\documentclass[11pt]{article}

\usepackage{amsmath}\usepackage{amsfonts}\usepackage{amssymb}\usepackage{graphicx}
\usepackage{ifthen}
\usepackage{lscape}
\usepackage{subfigure}
\usepackage{setspace}
\usepackage[left=1in,top=1.5in,right=1in,bottom=1in]{geometry}
\usepackage{verbatim}
\usepackage{mdwlist}
\usepackage{xcolor}
\usepackage{hyperref}

\usepackage{amsthm}
\newtheorem{theorem}{Theorem}

\renewcommand{\paragraph}[1]{}

\date{}
\author{Ping Lin\thanks{Division of Mathematics,
        University of Dundee, 23 Perth Road, Dundee, Scotland DD1 4HN, UK
        ({\tt plin@maths.dundee.ac.uk}).}
        \and
        Alexander V. Shapeev\thanks{Section of Mathematics, Swiss Federal Institute of Technology (EPFL), CH-1015, Lausanne, Switzerland
        ({\tt alexander.shapeev@epfl.ch}).}
        }

\title{Energy-Based Ghost Force Removing Techniques for the Quasicontinuum Method}

\newcommand{\Etot}{E_{\rm tot}}
\newcommand{\Eext}{E_{\rm ext}}
\newcommand{\Il}{\mathcal{I}_{\rm l}}
\newcommand{\In}{\mathcal{I}_{\rm n}}
\newcommand{\Inod}{\mathcal{I}_{\rm nod}}
\newcommand{\uQC}{\boldsymbol{u}^{\rm QC}}

\begin{document}
\sloppy

\maketitle

%\section*{TODO}
%
%\begin{itemize}
%\item alignment of formulae
%\end{itemize}

\begin{abstract}
This paper studies numerical methods for accurate treatment of the interface between the local and the nonlocal region in a QC approximation of atomistic materials.
Only the energy-based methods are considered.
Particularly, a quasicontinuum projection (QCP) method based on the idea of finite elements is shown to be accurate and efficient for this problem.
We analyse the QCP method and study its relation to the existing methods, such as the quasinonlocal quasicontinuum method and the geometrically consistent reconstruction-based method.
The analysis and the results of numerical tests confirm that the projection-based QC method successfully removes the ghost force with the same computational cost as the other methods.
In all computed examples the error of QCP is either the same or lower as the error of the other methods.
The performance of these methods in treating interfaces of elements in the local region is also examined.
\end{abstract}

\pagestyle{myheadings}
\thispagestyle{plain}
\markboth{PING LIN AND ALEXANDER V. SHAPEEV}{GHOST FORCE REMOVING FOR QC}

%\section*{Notations}
%
%\begin{description}
%\item[$N$] total number of atoms
%\item[$u_i$] positions of atoms
%\item[$f_i$] external force acting on atoms
%\item[$\Eext$] potential energy of the external force $\boldsymbol{f}$
%\item[$\Etot$] energy of atomistic interaction of the system
%\item[$\varphi$] interatomic interaction potential
%\item[$n$] interatomic interaction distance ($n=1$ is nearest neighbor interaction, $n=2$ is the next nearest neighbor interaction, etc).
%\item[$\Pi$] total potential energy of the atomistic material
%\item[$\Pi^{\rm XXX} / \Eext^{\rm XXX} / \Etot^{\rm XXX}$] potential energies associated with method ``XXX''.
%\item[$\epsilon$] characteristic interatomic distance
%\item[$\uQC$] solution by a quasicontinuum approximation
%\item[$\Il$] indices of ``local'' atoms (i.e.\ atoms subjected to smooth deformation)
%\item[$\In$] indices of ``nonlocal'' atoms (i.e.\ atoms subjected to non-smooth deformation)
%\item[$\Inod$] indices of nodal atoms (i.e.\ nodes of the triangulation of the domain)
%\item[${[v_i]}_i$] a vector whose components are $v_i$
%\item[${[m_{ij}]}_{ij}$] a matrix whose components are $m_{ij}$
%\end{description}

\section{Introduction}

\paragraph{ANOTHER VERSION}
Many fundamental problems in science and engineering can be modeled by partial differential equations, which rely on the macroscopic or continuum description of the problem.
In spite of tremendous success of continuum models, they also have certain limitations.
For example, when modeling macroscopic materials with microscopic defects, a continuum model will be less accurate near the defect and may give even wrong predictions, because the continuum model neglects the microscopic features of the defect.
Therefore, one might be tempted to switch to a complete atomistic model that takes the full atomistic description into account.
However, this is not an efficient strategy of modeling macroscopic materials not only because large atomistic systems are too large to handle even with the existing most powerful computers, but also because the solutions we obtain will likely contain too much information that is of little interest.
Thus, a combination of both atomistic and continuum descriptions is required for successful modeling of materials with defects.
Such an approach is called atomistic/continuum multiscale modeling.

Atomistic/continuum models adopt the point of view that there is an underlying atomistic model of the material which is the ``correct'' description or exact solution of the material problem.
This could be a quantum-mechanically based description such as density functional theory, but in practice it has been based primarily on classical-mechanics atomistic models with semi-empirical interatomic potentials.
The problem thus has the two length scales: the length of the whole material and the typical interatomic distance.
The ratio between these scales, $\epsilon$, is usually much less than $1$.
The bulk behavior of the material (at least in the case of small deformations) can be described by PDEs of continuum mechanics.
However, if one needs to model the material defects of the size $O(\epsilon)$, then one may require employing a fully atomistic model near those defects.

A representative of atomistic/continuum models is the quasi-continuum approximation (QC) which is becoming a popular multiscale technique for simulating static properties of polycrystalline materials.
It was put forward by Tadmor, Ortiz and
Phillips \cite{TOP} and later was updated in a review paper by
Miller and Tadmor \cite{MT}.
The idea of this approach is that we consider at the macroscopic scale the region in the material (called the local approximation region) where no serious defects occur and the theory of continuum material elasticity (or a coarser mesh approximation) may apply.
In the (nonlocal) region where serious defects occur, much finer mesh or atomistic model has to be used.
The QC method is usually described in a finite element framework.

While a significant body of knowledge about QC and related models and their numerical successes has been accumulated, not much has been reported on the analysis of these models.
The convergence rate of the QC method has been tested in \cite{KO2001} via computational results of a specific nano-indentation problem.
Other numerical demonstrations can also be found in the literature, e.g.\ \cite{SMTRPO1999,TOP}.
Rigorous numerical convergence analysis of the QC method has been published in \cite{Lin2003} for a one-dimensional case in the absence of external forces.
Later studies with certain external forces have been conducted in \cite{DLuskin,EM1,EM2,Ortner2008} for one-dimensional and occasionally two-dimensional cases.
An analysis of the method for two-dimensional cases under a few physically reasonable assumptions for solid materials has been given in \cite{linsinum}.
Other recent developments on the QC method and other atomistic/continuum models from the mathematical point of view include analysis in \cite{EM2} concerning the Cauchy-Born rule and in the context of the heterogeneous multiscale method \cite{EE}, in \cite{BLL} for excluding a spurious finite element effect in a prototypical multiscale model, in \cite{E2006} for a necessary and sufficient condition of the uniform first order accuracy, in \cite{ALuskin,OPB2005} for desirable a posteriori error estimators associated with the QC method, in \cite{LuskinOrtner08} for analysis of the cluster-summation rule,
%in \cite{GunzburgerZhang} ({\bf to be put in .bib file later}) for higher order finite elements,
and in \cite{BPBGL} for blending techniques in the mixed atomistic and continuum domain.

A problem having being associated with the QC method since its invention is the error at interfaces (edges or face) of elements and at the interface of local and nonlocal regions.
In \cite{linsinum} the error at edges of elements has been estimated in local regions for QC taking representative atoms at the center of finite elements.
Such an error increases with the number of edges, so at some critical point a mesh refinement would increase the numerical error even in the local QC case.
The error at the interfaces of local/nonlocal regions may be more difficult to study since it is usually associated with material defects.
Effort has been focused on unphysical forces (so-called ghost forces) at the local/nonlocal interfaces.
Such forces exist in the original QC method, unless only nearest neighbor interactions are taken at the local/nonlocal interfaces.
A correction force has been introduced in \cite{SMTRPO1999} and the method was called ``force-based correction''.
Convergence of an iterative procedure based on this version of QC has been shown in \cite{DLuskin}.
Despite its ability to eliminate the ghost force, the method is not associated with the ``correct'' potential energy and thus leads to an issue with energy conservation and possibly large numerical error (see e.g.\ \cite{Shimokawa2004}).
This is why we focus on the energy-based approach in this paper.
%But it does not mean that the force-based approach is not working.
The reader can refer to \cite{Dobson2009b-preprint,Ming08} for recent analysis of the force-based approach.
A seamless coupling in the interface region has been sought in \cite{Shimokawa2004} based on the changing the interaction energy of certain atoms and the resulting method is known as quasi-nonlocal method (QNL).
A convergence analysis for the linearized one-dimensional periodic case has been given in \cite{Dobson2009a-preprint} with the assumption of second nearest neighbor interaction.
In general the QNL method does remove the ghost force but is limited to potentials with a relatively short range.
In \cite{E2006,EMY08} it was shown that QNL works up to the second nearest neighbor interactions.
In \cite{EMY08} the accuracy of a more sophisticated approach which is based on what was called a geometrically consistent reconstruction (GCR) introduced in \cite{E2006} was studied in detail.
It was also claimed there that the GCR-based method includes QNL as a special case.
All these methods consider only reducing errors at the local-nonlocal interface.
However, the problem of reducing the error at edges of elements in the local region has not been addressed.

In this paper we aim to reexamine these existing methods, particularly focusing on existing versions of the QC method which is not originally proposed to remove the ghost force.
The method is based purely a finite element idea and often appears in literature (e.g.\ \cite{linplechac,Ortner2008,rudd2,RB2000}).
It is also briefly mentioned in \cite{RB05} for a linearized model that the method does not have the ghost force.
We shall call it the quasicontinuum projection method (QCP) in this paper.
Under a one-dimensional setting of the problem we will carefully study relationships between QNL, GCR and QCP, with a particular attention on the way these methods remove the ghost force.
We will show that QCP is not equivalent to QNL, not equivalent to the particular instance of GCR presented in \cite{E2006}, but there exist particular instances of GCR which are equivalent to QCP.
In addition, the QCP method offers a natural way of accurate treatment of edges of elements in the local region and thus eliminates the error at element interfaces.
This makes it an essentially more accurate method than the other ghost-force removing methods.
Based on that we conclude that QCP is the easiest and the most efficient way to remove the ghost force and the element edge errors in QC.
Furthermore, in favorable cases it does not introduce more computational costs in comparison with QNL and GCR.
In unfavorable cases, for instance when a 2D or 3D atomistic lattice is not aligned with the QC triangulation, QCP can be modified in the local region (while preserving the same accuracy in the local-nonlocal interface) so that the number of operation is of the same order as QNL or GCR. 
%We hope that this study would bring the understanding of this important issue of the QC method to a new level.

The paper is organized as follows.
The problem of an equilibrium of an atomistic material is introduced in section \ref{sec:problem-formulation} under the two settings: the general setting and the one-dimensional periodic setting.
It will be shown that QCP has no ghost force under the general setting in subsection \ref{sec:QCP-and-ghost-force}.
QC, QNL and GCR, their relations to each other, and their treatment of the ghost force are examined under the periodic setting in section \ref{sec:QC}.
QCP is then introduced and compared with QNL and GCR in subsections \ref{sec:QCP:general}--\ref{sec:QCP-implementation}.
Numerical experiments under the both one-dimensional periodic setting and general setting are given in section \ref{sec:results} to show the performance of these methods.
The results are discussed and summarized in section \ref{sec:discussion-conclusion}.

\section{Problem Formulation}\label{sec:problem-formulation}

\paragraph{intro}
%???link with the Introduction
We first describe the general problem formulation and then will consider the 1D period setting (as considered, for instance, in \cite{Dobson2009a-preprint}).

\subsection{General Setting}

\paragraph{introduction of atoms and atomistic energy}
Consider a problem of finding an equilibrium configuration of atoms in an atomistic material in space $\mathbb{R}^d$.
Denote the degrees of freedom (which are chosen to be the coordinates of atoms) by a vector
\[
\boldsymbol{u} = \{u_i\} , \quad (1\le i\le N).
\]
To model the equilibrium of atoms, we introduce the potential energy of the atomistic system
\begin{equation}
\Pi(\boldsymbol{u}) = \Etot(\boldsymbol{u}) + \Eext(\boldsymbol{u}),
\label{eq:general_Pi}
\end{equation}
which is a sum of the internal energy of atoms $\Etot(\boldsymbol{u})$ and the external potential $\Eext(\boldsymbol{u})$, the later will be considered to be a linear functional of $\boldsymbol{u}$ (which means that the external forces on atoms do not depend on their positions).

\paragraph{equilibrium}
In the above notations the problem of finding the equilibrium configuration of atoms is written as
\begin{equation}
\frac{\partial \Pi}{\partial u_i} = 0 \quad (i=1,2,\ldots,N).
\label{eq:original_problem}
\end{equation}
Sometimes, instead of simply an equilibrium configuration, one specifically requires a stable equilibrium configuration, which is formally written as
\begin{equation}
\textnormal{$\boldsymbol{u}$ is a local minimizer of $\Pi(\boldsymbol{u})$}.
\label{eq:original_stable_problem}
\end{equation}

The equilibrium equations \eqref{eq:original_problem} are often solved by Newton's iterative method which in matrix form reads
\[
\left[
	\frac{\partial^2 \Pi(\boldsymbol{u}^{n})}{\partial u_i \partial u_j}
\right]_{ij}
(\boldsymbol{u}^{n+1}-\boldsymbol{u}^{n}) = - \left[ \frac{\partial \Pi(\boldsymbol{u}^{n})}{\partial u_i} \right]_i.
\]
To assemble this system of equations one should effectively compute the right-hand side vector components $\frac{\partial \Pi}{\partial u_i}$ and the stiffness matrix components $\frac{\partial^2 \Pi}{\partial u_i \partial u_j}$.
The algorithm of computing these would depend on a particular problem setting.
Below we elaborate the 1D periodic setting, since we find it to be the easiest formulation for analysis and comparison of different QC approaches, especially if one's focus is the ghost force.

\subsection{1D Periodic Setting}

\paragraph{intro}
In this subsection we describe the problem formulation of finding an equilibrium of an atomistic material in the 1D periodic setting.
We consider the periodic boundary conditions, as it is done in \cite{Dobson2009a-preprint}.
This is done to avoid difficulties arising from presence of the boundary of the atomistic material.
Otherwise, the boundary of an atomistic material, unless properly treated, would contribute an additional error to the numerical solution.
%, since the focus of the present work is on elimination of the error on the interface between local and nonlocal regions (hereinafter called ``{\itshape the interface error}'').
%In subsection \ref{sec:QCP-implementation} we will mention implementation of the projection method in the general case and in subsection \ref{sec:2dtest} we will present the numerical results in the 2D non-periodic case.

\paragraph{introduction of atoms and atomistic energy}
Consider an atomistic material in one dimension with atom positions $u_i$, $(-\infty<i<\infty)$.
We assume that the material behaves periodically with a period of length one over $N$ atoms:
\[
u_{i+N} = u_i+1 \quad (-\infty<i<\infty).
\]
Thus, the degrees of freedom of such system are described by
\[
\boldsymbol{u} = \{u_i\}, \quad (1\le i\le N).
\]
The energy of atomistic interaction of the system (summed for the atoms over one period) is then
\begin{equation}
\Etot(\boldsymbol{u}) =
%\epsilon\sum_{i=1}^{N-1} \sum_{j=-\infty}^{\infty} \varphi\left(\frac{u_i-u_j}{\epsilon}\right) =
\frac{\epsilon}{2} \sum_{i=1}^{N} \sum_{\substack{j=-\infty \\ j\ne i}}^{\infty} \varphi\left(\frac{u_i-u_j}{\epsilon}\right),
\label{eq:Eint}
\end{equation}
where $\varphi(r/\epsilon)$ is the potential of interaction of atoms at distance $r$.
For compactness of notations, we define the potential $\varphi(z)$ for negative $z$ as $\varphi(z) = \varphi(-z)$.
We assume that the potential $\varphi(z)$ vanishes for $|z|$ large enough, so that we need to consider at most $n$ neighboring atoms in the interaction energy:
\[
\Etot(\boldsymbol{u}) =
\frac{\epsilon}{2} \sum_{i=1}^{N} \sum_{\substack{j=i-n\\ j\ne i}}^{i+n} \varphi\left(\frac{u_i-u_j}{\epsilon}\right).
\]

\paragraph{introduction of potential energy, problem formulation}
The potential energy of the external force $\boldsymbol{f}$ is
\begin{equation}
\Eext(\boldsymbol{u}) = -\epsilon\sum_{i=1}^{N} f_i u_i.
\label{eq:Eext}
\end{equation}
The forces $f_i$ on each atom are given and considered to be independent of atom positions $u_i$.
For the problem to be well-posed, the sum of all forces is assumed to be zero:
\[
\sum_{i=1}^{N} f_i = 0.
\]
If the sum is not zero, then the total energy will not be bounded from below due to the periodicity.

\section{QC approximation}\label{sec:QC}

\paragraph{introduce QC approximation: nodal atoms}
Usually, the characteristic distances on which the external force $\boldsymbol{f}$ varies are much larger than the characteristic interatomic distance $\epsilon$.
In this case we do not need to model all degrees of freedom of the system associated with each individual atom and can employ a QC approximation.
It consists in choosing the nodal atoms $\Inod = \{i_1,i_2,\ldots,i_K\}$ and regarding the positions of the non-nodal atoms to be known via a linear interpolation:
\begin{equation}
u_i = \frac{i_k-i}{i_k-i_{k-1}} u_{i_{k-1}} + \frac{i-i_{k-1}}{i_k-i_{k-1}} u_{i_k}
\quad (i_{k-1}< i< i_k),
\quad k=1,\ldots,K.
\label{eq:QC_constraints}
\end{equation}
Here we defined $i_0 = i_K - N$ by periodicity.
Then all distances between neighboring atoms positioned between $i_{k-1}$ and $i_k$ are equal to the same value
\[
u_i - u_{i-1} = \frac{u_{i_k}-u_{i_{k-1}}}{i_k-i_{k-1}}
\quad (i_{k-1}<i\le i_k).
\]
Thus, the unknowns in a QC approximation are
\begin{equation}
\uQC = [ u_{i_1}, u_{i_2}, \ldots, u_{i_K}  ] = \left[u_{i_k}\right]_k.
\label{eq:uQC}
\end{equation}

\paragraph{original QC: energy approximation}
After having defined the nodal atoms, we need to approximate the original problem \eqref{eq:original_problem}.
The original QC method (see e.g.\ \cite{linsinum,MT,TOP}) is based on the assumption that the energy of the segment $i_{k-1}\le i\le i_k$ can be approximated by the continuum energy via the so-called Cauchy-Born rule:
\[
\frac{\epsilon}{2} \sum_{i=i_{k-1}}^{i_k} \sum_{\substack{j=-n\\ j\ne i}}^{n} \varphi\left(\frac{u_i-u_j}{\epsilon}\right)
\approx
\epsilon(i_k-i_{k-1})\,\Phi\left(\frac{u_{i_k}-u_{i_{k-1}}}{\epsilon(i_k-i_{k-1})}\right)
\]
where the interatomic energy density $\Phi$ is defined as
\begin{equation}
\Phi(z) = \frac{1}{2} \sum_{\substack{m=-n \\ m\ne 0}}^{n} \varphi\left(m z\right)
 = \sum_{m=1}^{n} \varphi\left(m z\right).
\label{eq:Phi}
\end{equation}
%It is based on the assumption that the number of atoms between $i_k$ and $i_{k-1}$ is large enough to neglect the fact that ...?

\paragraph{external force is localized}
In the present paper, we are focused on modeling a localized irregularity of the atomistic material, which requires keeping track of individual atoms in the neighborhood of the irregularity (called nonlocal QC) while using a local QC approximation elsewhere.
To model such situation in 1D, we consider the external force $\boldsymbol{f}$ to be zero everywhere except for the atoms $M,M+1,\ldots,M+P-1$:
\[
f_i = 0 \quad (i=1,\ldots, M-1, M+P,\ldots,N),
\]
where $M$ is chosen to be $M=\lfloor N/2 \rfloor$.
%We note that this is a simplified assumption.
%A more general assumption is that $\boldsymbol{f}$ varies smoothly everywhere except a limited number of atoms.

\begin{figure}
\begin{center}
\includegraphics{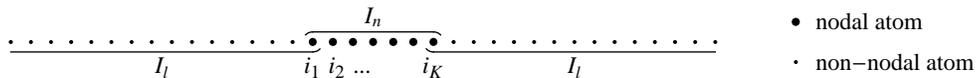}
\caption{Illustration of local and nonlocal regions}
\label{fig:illustration_regions}
\end{center}
\end{figure}

\begin{figure}
\begin{center}
\includegraphics{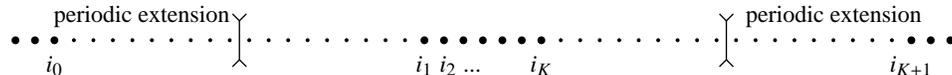}
\caption{Illustration of periodic extension}
\label{fig:illustration_periodic}
\end{center}
\end{figure}

\paragraph{introduction of local and nonlocal regions}
A common approach to this problem is the following \cite{Shenoy1999}.
Introduce the nonlocal region
\[
\In=\{i: \ M-L\le i\le M+P-1+L\}
\]
and the local region
\[
\Il=\{1,2,\ldots,N\}\setminus \In.
\]
The nonlocal region $\In$ is chosen as an extension of the region where $f_i\ne 0$ by $L$ atoms in each direction.
The nodal atoms $\Inod = \{i_1,i_2,\ldots,i_K\}$ are chosen so that $\In \subset \Inod$.
Since we assumed the external force $\boldsymbol{f}$ to be equal to zero in the local region $\Il$, the deformations will be well-defined by a linear function and therefore nodal atoms are not required to be added inside $\Il$ (see figs.\ \ref{fig:illustration_regions} and  \ref{fig:illustration_periodic}).
We stress again that this specific problem setting is chosen to study the relationships between ghost-force removing strategies.
In a more general setting, local and nonlocal regions, as well as the mesh refinement in local regions, may be determined through a certain adaptive strategy (see e.g.\ \cite{MT,TOP}).

\paragraph{formulae of the QCE}
Thus, the expression of the total internal energy of the standard QC method is
\[
\Etot^{\rm QCE} = E_{\rm local}^{\rm QCE} + E_{\rm nonlocal}^{\rm QCE},
\]
where
\[
E_{\rm local}^{\rm QCE} = \epsilon(i_1-i_0)\,\Phi\left(\frac{u_{i_1}-u_{i_0}}{\epsilon(i_1-i_0)}\right)
\]
for $\Phi$ defined by \eqref{eq:Phi}, and
\[
E_{\rm nonlocal}^{\rm QCE} = \frac{\epsilon}{2}\sum\limits_{i=i_1}^{i_K}\sum\limits_{\substack{j=i-n \\ j\ne i}}^{i+n} \varphi\left(\frac{u_i-u_j}{\epsilon}\right).
\]
The respective system of equations to be solved is
\begin{equation}
\begin{array}{l} \displaystyle
\frac{\partial \Pi^{\rm QCE}}{\partial u_{i_k}} =
\frac{\partial \Etot^{\rm QCE}}{\partial u_{i_k}} + \frac{\partial \Eext}{\partial u_{i_k}} = 0 \quad (k=1,2,\ldots,K),
\end{array}
\label{eq:QCE_equations}
\end{equation}
which is the shortened version of the original system \eqref{eq:original_problem}.

\subsection{Ghost Force, QNL Method, and GCR method}

\paragraph{ghost force}
The traditional QC method for atomistic material modeling, though seems to be very intuitive, has an essential drawback: the accuracy of approximation in the neighborhood of the interface between $\In$ and $\Il$ falls down to zeroth order.
This phenomenon of loss of accuracy on the interface between the regions is commonly interpreted as the ghost force (see e.g.\ \cite{CurtinMiller03,MT,Shenoy1999}).

\paragraph{ghost force illustration}
The nature of the ghost force can be illustrated in the case of second nearest neighbor interaction $n=2$.
Consider the equidistant lattice $u_i = \epsilon z i$ with the zero external force $\boldsymbol{f}=0$ (and hence zero $\Eext$).
In this case the force on each atom should be zero (due to reflection symmetry).
Compute the force on the atom $i_2$:
\[
\frac{\partial \Etot^{\rm QCE}}{\partial u_{i_2}} =
\frac{\partial E_{\rm nonlocal}^{\rm QCE}}{\partial u_{i_2}} + \frac{\partial E_{\rm local}^{\rm QCE}}{\partial u_{i_2}}
= \frac{\epsilon}{2} \frac{\partial}{\partial u_{i_2}}\sum\limits_{i=i_1}^{i_K}\sum\limits_{\substack{j=i-n \\ j\ne i}}^{i+n} \varphi\left(\frac{u_i-u_j}{\epsilon}\right) + 0
\]
\[
=
\frac{\epsilon}{2}\left(-\varphi'\left(\frac{u_{i_1}-u_{i_2}}{\epsilon}\right)
+ \varphi'\left(\frac{u_{i_2}-u_{i_2-2}}{\epsilon}\right)+\varphi'\left(\frac{u_{i_2}-u_{i_1}}{\epsilon}\right)+\varphi'\left(\frac{u_{i_2}-u_{i_3}}{\epsilon}\right)
\right.
\]
\[
\left.
+\varphi'\left(\frac{u_{i_2}-u_{i_4}}{\epsilon}\right)
-\varphi'\left(\frac{u_{i_3}-u_{i_2}}{\epsilon}\right)-\varphi'\left(\frac{u_{i_4}-u_{i_2}}{\epsilon}\right)\right),
\]
which upon substitution $u_i = \epsilon z i$ yields
\[
\frac{\partial \Etot^{\rm QCE}}{\partial u_{i_2}} = \frac{\epsilon}{2} \varphi'\left(2 z\right).
\]
The reason for not getting zero force is that we are missing the term
\[
\frac{\epsilon}{2} \varphi\left((u_{i_2-2}-u_{i_2})/\epsilon\right)
\]
in the expression for $\Etot^{\rm QCE}$ (specifically in $E_{\rm local}^{\rm QCE}$).

\paragraph{QNL method}
There are several methods to eliminate the ghost force, most of which can be divided into the two categories: force-based methods (which modify the right-hand side of the equilibrium equations \eqref{eq:QCE_equations}) and energy-based methods, which modify the interaction potential $\Etot^{\rm QCE}$.
As mentioned in Introduction, we will not discuss any force-based methods but focus only on the energy-based methods.
One of the energy-based methods proposed to eliminate the ghost force was the so-called quasinonlocal quasicontinuum method, or QNL in short \cite{Shimokawa2004}.
The idea of the QNL method is to modify the interaction of atom $i_2$ (and the respective atom on the other side of the nonlocal region $i_{K-1}$) in order to eliminate the ghost force.
Namely, instead of the term $\varphi\left(\frac{u_{i_2}-u_{i_2-2}}{\epsilon}\right)$ we introduce the term $\varphi\left(\frac{2 (u_{i_2}-u_{i_2-1})}{\epsilon}\right)$, just as in the Cauchy-Born extrapolation, but for the nonlocal atom $i_2$.
In this case the atom $i_2$ is called a quasinonlocal atom, since though it formally lies in the nonlocal region, it interacts with atom $i_2-2$ like a local atom.
Thus, the expression of the energy of the QNL method is
\[
\Etot^{\rm QNL} = E_{\rm local}^{\rm QNL} + E_{\rm nonlocal}^{\rm QNL} + E_{\rm quasinonlocal}^{\rm QNL},
\]
where
\[
E_{\rm local}^{\rm QNL} = \epsilon(i_1-i_0)\,\Phi\left(\frac{u_{i_1}-u_{i_0}}{\epsilon(i_1-i_0)}\right)
\]
for $\Phi$ defined by \eqref{eq:Phi}, and
\[
E_{\rm nonlocal}^{\rm QNL} = \frac{\epsilon}{2}\sum\limits_{\substack{i=i_1 \\ i\notin \{i_2,i_{K-1}\}}}^{i_K}\sum\limits_{\substack{j=i-n \\ j\ne i}}^{i+n} \varphi\left(\frac{u_i-u_j}{\epsilon}\right),
\]
and finally
\[
E_{\rm quasinonlocal}^{\rm QNL} = \frac{\epsilon}{2}
\left(
	\varphi\left(\frac{2(u_{i_2}-u_{i_2-1)}}{\epsilon}\right) +
	\varphi\left(\frac{u_{i_2}-u_{i_2-1}}{\epsilon}\right)
\right.
\]
\[
\left.
+
	\varphi\left(\frac{u_{i_2}-u_{i_2+1}}{\epsilon}\right) +
	\varphi\left(\frac{u_{i_2}-u_{i_2+2}}{\epsilon}\right)
\right)
\]
\[
+
\frac{\epsilon}{2}
\left(
	\varphi\left(\frac{u_{i_{K-1}}-u_{i_{K-1}-2}}{\epsilon}\right) +
	\varphi\left(\frac{u_{i_{K-1}}-u_{i_{K-1}-1}}{\epsilon}\right)
\right.
\]
\[
\left.
+
	\varphi\left(\frac{u_{i_{K-1}}-u_{i_{K-1}+1}}{\epsilon}\right) +
	\varphi\left(\frac{2 (u_{i_{K-1}}-u_{i_{K-1}+1})}{\epsilon}\right)
\right).
\]
%One can check that for the QNL method there is no ghost force for second nearest neighbor interaction $n=2$.
%However in the case $n>2$ there is the ghost force.

\paragraph{ghost force for QNL}
The QNL method has been designed to have no ghost force for the case $n=2$, which means that the force exerted on the atoms in the equispaced lattice $u_i = \epsilon z i$ is exactly zero.
One, however, can do a more detailed analysis: consider an arbitrary lattice satisfying the QC constraints \eqref{eq:QC_constraints} and compare the QNL energy with the exact atomistic energy:
\[
\Etot^{\rm QNL} - \Etot =
\frac{\epsilon}{2} \varphi\left(\frac{2 (u_{i_1+1}-u_{i_1})}{\epsilon}\right)
+
\frac{\epsilon}{2} \varphi\left(\frac{2 (u_{i_1}-u_{i_1-1})}{\epsilon}\right)
-\epsilon \varphi\left(\frac{u_{i_1+1}-u_{i_1-1}}{\epsilon}\right)
\]
\[
+ \textnormal{ respective terms for atoms around $i_K$}.
\]
Defining $D^2u_i = \frac{u_{i+1}-2 u_i+u_{i-1}}{\epsilon^2}$ and expanding $\Etot^{\rm QNL}-\Etot$ in Taylor series w.r.t.\ $u_{i_2}$ around its extrapolated value $2 u_{i_1}-u_{i_1-1}$ (the respective terms for atoms around $i_K$ are expanded w.r.t.\ $u_{i_K-1}$) yields
\[
\Etot^{\rm QNL} - \Etot =
\frac{\epsilon}{2} \left(\epsilon\,\frac{u_{i_1+1}-2 u_{i_1}+u_{i_1-1}}{\epsilon^2}\right)^2 \varphi''\left(\frac{2 (u_{i_1}-u_{i_1-1})}{\epsilon}\right)
\]
\[
+ \epsilon O\left(\epsilon\,\frac{u_{i_1+1}-2 u_{i_1}+u_{i_1-1}}{\epsilon^2}\right)^4
+ \textnormal{ respective terms for atoms around $i_K$}
\]
\[
=
\frac{\epsilon}{2} \left(\epsilon D^2u_{i_1}\right)^2 \varphi''\left(\frac{2 (u_{i_1}-u_{i_1-1})}{\epsilon}\right)
+ \epsilon O\left(\epsilon D^2u_{i_1}\right)^4
\] \[
+
\frac{\epsilon}{2} \left(\epsilon D^2u_{i_K}\right)^2 \varphi''\left(\frac{2 (u_{i_K}-u_{i_K+1})}{\epsilon}\right)
+ \epsilon O\left(\epsilon D^2u_{i_K}\right)^4.
\]
\paragraph{difference is small near the interface}
One can see that if the interface atom $i_1$ (as well as $i_K$) lies in the region where the solution is smooth (i.e.\ $D^2 u_{i_1}$ is bounded, then the difference $\Etot^{\rm QNL}-\Etot$ is small.
Moreover, the force on the atom $i_2$ is now
\[
\frac{\partial \Etot^{\rm QNL}}{\partial u_{i_2}}
-
\frac{\partial \Etot}{\partial u_{i_2}}
=
\frac{\epsilon}{2} \left(D^2u_{i_1}\right)^2 + O\left(\epsilon^3 (D^2u_{i_1})^3\right),
\]
which is $O(\epsilon)$ if the solution is smooth near the interface.
In the next section we will show that the QCP method has exactly the energy $\Etot$, which means that it is not equivalent to the QNL method.

\paragraph{ELY's method}
E, Lu, and Yang \cite{E2006} proposed a way to remove the ghost force for an arbitrary interaction distance $n$, based on what they call a geometrically consistent reconstruction (GCR).
They interpreted the change of the interaction of certain atoms as changing the reconstruction of atoms: for instance, in the QNL method, in the interaction of atom $i_2$ with $i_2-2$, the position of atom $i_2-2$ is reconstructed as $u_{i_2}+2 (u_{i_1}-u_{i_2})$.
The idea in \cite{E2006} was to change the interaction (i.e.\ reconstruction) of more than one atom near the interface (the number of atoms whose reconstruction is changed depends on $n$) so that the ghost force does not appear for a given $n$.
The reconstruction was sought by the indeterminate coefficients method in the form of a linear combination of the nonlocal reconstruction and the Cauchy-Born reconstruction:
\begin{equation}
u^{\rm reconstructed}_j = C^{\rm GCR}_{ij} u_j + \left(1-C^{\rm GCR}_{ij}\right) (u_i + (j-i) (u_{i+{\rm sgn}(j-i)}-u_i)).
\label{eq:GCR-reconstruction}
\end{equation}
The requirement of such reconstruction to have no ghost force results in a class of methods.
E, Lu, and Yang gave particular instances of such methods in \cite{E2006}.
For $n\le 3$ the coefficients $C_{ij}$ given in table I in \cite[p.\ 214115-8]{E2006} are:
\[
C^{\rm GCR}_{ij} =
\left\{
\begin{array}{lll}
1 & & \textnormal{$(i,j)$ = $(i_1-1,i_1)$ or $(i_K+1,i_K)$}, \\
1 & & \textnormal{$(i,j)$ = $(i_1-1,i_1+1)$ or $(i_K+1,i_K-1)$}, \\
1 & & \textnormal{$(i,j)$ = $(i_1-1,i_1+2)$ or $(i_K+1,i_K-2)$}, \\
2/3 & & \textnormal{$(i,j)$ = $(i_1-2,i_1+1)$ or $(i_K+2,i_K-1)$}, \\
1/3 & & \textnormal{$(i,j)$ = $(i_1+1,i_1-2)$ or $(i_K-1,i_K+2)$}, \\
C^{\rm QCE}_{ij} & & \textnormal{otherwise},
\end{array}
\right.
\]
where the coefficients $C^{\rm QCE}_{ij}$ are defined as follows:
\[
C^{\rm QCE}_{ij} =
\left\{
\begin{array}{lll}
1 & & \textnormal{$i<i_1$ or $i>i_K$}, \\
0 & & \textnormal{$i_1<i<i_K$}, \\
1 & & \textnormal{$i\in \{i_1,i_K\}$ and $i_1<j<i_K$}, \\
0 & & \textnormal{$i\in \{i_1,i_K\}$ and ($i<i_1$ or $i>i_K$)}.
\end{array}
\right.
\]

\paragraph{example for $n=2$ of ELY's method}
For example, as is implied from table I in \cite[p.\ 214115-8]{E2006}, for $n=2$
%(the case which we will use the method of \cite{E2006} in the present numerical experiments,
the atomistic energy (for brevity written below as a correction to $\Etot^{\rm QCE}$) is
\[
\Etot^{\rm GCR} = \Etot^{\rm QCE}
+
\left(
	\varphi\left(\frac{u_{i_1-1}-u_{i_1+1}}{\epsilon}\right)
-
	\varphi\left(\frac{2 (u_{i_1-1}-u_{i_1})}{\epsilon}\right)
\right)
\]
\[
+
\left(
	\varphi\left(\frac{u_{i_K+1}-u_{i_K-1}}{\epsilon}\right)
-
	\varphi\left(\frac{2 (u_{i_K+1}-u_{i_K})}{\epsilon}\right)
\right).
\]
Note that here we ignored another correction term
\[
\varphi\left(\frac{u_{i_1}-u_{i_1-2}}{\epsilon}\right)-\varphi\left(\frac{2 (u_{i_1}-u_{i_1-1})}{\epsilon}\right)
\]
 which equals zero identically under the QC constraints \eqref{eq:QC_constraints}, but is implied in the table in \cite{E2006}.
It is worth noting that this particular instance of the GCR method is not equivalent to the QNL method for $n=2$ (although \cite{E2006} claimed it is equivalent): the QNL does not change the interaction of the local atom $i_1-1$, whereas that particular GCR method does.
However, it is true that the QNL method is a particular instance of the GCR method.
Particularly, if the one changes the coefficients $C^{\pm}_i$ to $C^{\pm}_{i+1}$ in table I of \cite{E2006} for $n=1$ and $n=2$, then one will obtain exactly the QNL method.
That is, in notations of this paper the QNL method's reconstruction is
\[
C^{\rm QNL}_{ij} =
\left\{
\begin{array}{lll}
0 & & \textnormal{$(i,j)$ = $(i_1,i_3)$ or $(i_{K-2},i_K)$}, \\
C^{\rm QCE}_{ij} & & \textnormal{otherwise},
\end{array}
\right.
\]
which is clearly different from $C^{\rm GCR}_{ij}$.

\paragraph{GCR energy is exact for $n=2$}
Another property of the instance of the GCR method in table I of \cite{E2006} which is worth noting is that atomistic energy of the GCR method identically equals to the exact atomistic energy: $\Etot^{\rm GCR} = \Etot$ for $n=2$ for a lattice satisfying \eqref{eq:QC_constraints}.
As we will see below, it is not the case for $n=3$, although one can find other instances of the GCR method, which would give the exact energy $\Etot$ on the QC lattice for a given $n$.

%\paragraph{example of ELY's method for $n=3$}
%For $n=3$ the presented GCP method has the following interaction energy
%\[
%\Etot^{\rm GCR} = \Etot^{\rm QCE} +
%\left(
%	\varphi\left(\frac{u_{i_1-1}-u_{i_1+1}}{\epsilon}\right)
%-
%	\varphi\left(\frac{2 (u_{i_1-1}-u_{i_1})}{\epsilon}\right)
%\right)
%\]
%\[
%-
%\epsilon\varphi\left(\frac{u_{i_1+1}-u_{i_1-2}}{\epsilon} \right)
%\]
%\[
%+
%\frac{\epsilon}{2}\varphi\left(\frac{5 (u_{i_2}-u_{i_1})}{3 \epsilon} + \frac{4 (u_{i_1}-u_{i_0})}{3 \epsilon(i_1-i_0)} \right)
%-
%\frac{\epsilon}{2}\varphi\left(\frac{u_{i_2}-u_{i_1}}{3 \epsilon} + \frac{8 (u_{i_1}-u_{i_0})}{3 \epsilon(i_1-i_0)} \right).
%\]
%
%%The first method proposed to eliminate the ghost force was the so-called quasinonlocal quasicontinuum method.
%%Below we describe
%%A number of papers suggest different approaches to eliminate the ghost force in order to get at least first order convergence of the QC approximation.

\paragraph{proposed approach}
GCR gives a general framework or criterion for removing the ghost force at the local-nonlocal interface.
Nevertheless, it involves a priori calculation of reconstruction parameters and the formulation of the method is different for any slight change of the problem, e.g.\ from $n=2$, $3$ to any given number.
In the present paper we thus seek and study a specific method simply based on the projection \eqref{eq:QC_constraints} requiring no a priori calculation of any parameters and meanwhile performing not worse than GCR in all occasions.
We will call it the QC projection method (QCP).
%focus more on treating the ghost force simply using the finite element approach based on the projection \eqref{eq:QC_constraints}, which we will call the QC projection method (QCP).

\section{QCP method}

\paragraph{intro}
The idea of projection was used by Rudd and Broughton in \cite{rudd2,RB2000,RB05} in their coarse-grained molecular dynamics method and the absence of the ghost force for such method was noticed by the same authors in \cite{RB05} for a linearized model.
A rather detailed analysis of such method can be found in \cite{Ortner2008} for a one-dimensional steady-state problem satisfying the Dirichlet boundary condition.
In this method, the full atomistic solution is ``projected'' onto the space of admissible deformations \eqref{eq:QC_constraints}, in the same sense as the solution to differential equations is projected onto the finite element space in Galerkin projection methods.
We therefore call it the quasicontinuum projection method (QCP).
%...propose a novel approach to treating the error on the interface between the nonlocal and the local region.
%We will give our own interpretation of this interface error by treating the QC approximation as the projection of the original problem \eqref{eq:original_problem} and interpreting the interface error as a conformity of the projection.
This approach can also be interpreted as introducing the constraints \eqref{eq:QC_constraints} to the original optimization problem \eqref{eq:original_problem}.
In this case the approximate solution $\uQC$ (cf.\ \eqref{eq:uQC}) will be the best possible QC solution, i.e.\ the solution which minimizes the energy $\Pi(\boldsymbol{u})$ among all $\boldsymbol{u}$ satisfying \eqref{eq:QC_constraints}.
%This is equivalent to the Ortner and S\"uli's version of the QC method \cite{Ortner2008} derived as a Galerkin projection of \eqref{eq:original_problem} on the subspace defined by \eqref{eq:QC_constraints}.

\paragraph{organization of the section}
In this section we first formulate the quasicontinuum projection (QCP) method in the general setting, then we formulate it in the 1D periodic setting for the convenience of discussing the relationship of QCP with the other ghost-force removing methods near the local-nonlocal interface.

\subsection{QCP in General Setting}\label{sec:QCP:general}

\paragraph{idea}
We want to approximately find the critical point of $\Pi(\boldsymbol{u})$ by reducing the number of degrees of freedom of $\boldsymbol{u}$ by considering the constraints \eqref{eq:QC_constraints}.
A natural approach would be to just restrict our search for the critical point to the subspace defined by the constraints \eqref{eq:QC_constraints}.

\paragraph{details}
Formally this can be written as follows.
Define the QC variables
\[
\uQC=[ u_{i_1}, u_{i_2}, \ldots, u_{i_K} ] = \left[u_{i_k}\right]_k,
\]
which are analogous to those in 1D periodic setting \eqref{eq:uQC}.
Define the QC reconstruction matrix $U^{\rm QC}$ in the following way:
\begin{equation}
\boldsymbol{u} = U^{\rm QC}\uQC,
\label{eq:matrix-U}
\end{equation}
which corresponds to the QC reconstruction \eqref{eq:QC_constraints}.
Then the QCP method can formally be written as
\begin{equation}
\frac{\partial}{\partial u_{i_k}}\,\Pi\left(U^{\rm QC}\uQC\right) = 0
\quad (1\le i\le K),
\label{eq:general_QC_equations}
\end{equation}
or alternatively, in relation to the stable equilibrium problem \eqref{eq:original_stable_problem} as
\[
\textnormal{$\uQC$ is a local minimizer of $\Pi\left(U^{\rm QC}\uQC\right)$}.
\]
Define the QCP approximation to the potential energy as
\begin{equation}
\Pi^{\rm QCP}\left(\uQC\right) = \Pi\left(U^{\rm QC}\uQC\right).
\label{eq:Pi-QCP}
\end{equation}

\paragraph{method}
The equilibrium QC equations \eqref{eq:general_QC_equations} can also be solved using the Newton's method
\begin{equation}
\left[
	\frac{\partial^2 \Pi\left(U^{\rm QC}\boldsymbol{u}^{{\rm QC},n}\right)}{\partial u_{i_k} \partial u_{i_l}}
\right]_{kl}
(\boldsymbol{u}^{{\rm QC},n+1}-\boldsymbol{u}^{{\rm QC},n}) = - \left[ \frac{\partial \Pi\left(U^{\rm QC}\boldsymbol{u}^{{\rm QC},n}\right)}{\partial u_k} \right]_k.
\label{eq:general_QCP_Newton}
\end{equation}
The implementation of the QCP method in the general setting will be discussed after further elaboration of the QCP method in the 1D periodic setting.

\subsection{QCP in 1D periodic setting}

\paragraph{intro}
In this subsection we formulate the QCP method in 1D periodic setting and analyze it.
%Nevertheless in subsection \ref{sec:QCP-implementation} we discuss the implementation of the QCP method and indicate how to extend this method to more general cases.
For that we take the specific form of the reconstruction operator $U^{\rm QC}$ which corresponds to formula \eqref{eq:QC_constraints}:
\[
(U^{\rm QC} \uQC)_i = \frac{i_k-i}{i_k-i_{k-1}} u_{i_{k-1}} + \frac{i-i_{k-1}}{i_k-i_{k-1}} u_{i_k}
\quad (i_{k-1}< i< i_k).
\]
Here we continue using the periodic extension of indices: $i_{K+1}=i_1$ and $i_0 = i_K$.
Notice that the matrix $U^{\rm QC}$ is sparse: each row of $U^{\rm QC}$ consists of at most two entries.

\paragraph{energy formula}
Now, according to \eqref{eq:general_QC_equations}, we need to write down the potential energy
\[
\Pi^{\rm QCP}(\uQC) = \Pi(U^{\rm QC} \uQC),
\]
which is nothing but the energy \eqref{eq:general_Pi} written in terms of positions of the nodal atoms $u_{i_1},\ldots,u_{i_K}$.
The external energy $\Eext$ is easy to deal with since $\boldsymbol{f}$ is zero for the non-nodal atoms.
To write down the atomistic interaction energy of QCP, we split it into the three parts, corresponding to interaction between nonlocal atoms ($i_1\le i\le i_K$), interaction between local atoms ($i<i_1$ or $i>i_K$), and interaction of nonlocal atoms with local atoms:
\begin{equation}
\Etot^{\rm QCP} = E_1^{\rm QCP} + E_2^{\rm QCP} + E_3^{\rm QCP}.
\label{eq:QCP-energy}
\end{equation}
The energy of interaction between nonlocal atoms $E_1^{\rm QCP}$ can easily be written down since all nonlocal atoms are nodal atoms:
\begin{equation}
E_1^{\rm QCP} =
\frac{\epsilon}{2} \sum_{\substack{i_1\le i,j\le i_K \\ 0<|i-j|\le n}} \varphi\left(\frac{u_i-u_j}{\epsilon}\right).
\label{eq:QCP-energy-1}
\end{equation}
The energy of interaction between local atoms $E_2^{\rm QCP}$ must be written in terms of the nodal atoms $u_{i_1}$ and $u_{i_K}$.
For this purpose, we denote $i_0 = i_K - N$ and will use the periodicity condition $u_{i-N} = u_i-1$ for writing down this energy.
We will also use $u_i-u_j = (i-j)\frac{u_{i_1}-u_{i_0}}{i_1-i_0}$.
Thus, the energy $E_2^{\rm QCP}$ is
\begin{equation}
\begin{array}{l} \displaystyle
E_2^{\rm QCP} =
\frac{\epsilon}{2} \sum_{\substack{i_0\le i,j\le i_1 \\ 0<|i-j|\le n}} \varphi\left(\frac{(i-j) (u_{i_1}-u_{i_0})}{\epsilon (i_1-i_0)}\right)
\\ \displaystyle
=
\frac{\epsilon}{2} \sum_{\substack{m=-n \\ m\ne 0}}^{n}(i_1-i_0+1-|m|) \varphi\left(\frac{m (u_{i_1}-u_{i_0})}{\epsilon(i_1-i_0)}\right)
\\ \displaystyle
=
\epsilon \sum_{m=1}^{n}(i_1-i_0+1-m) \varphi\left(\frac{m (u_{i_1}-u_{i_0})}{\epsilon(i_1-i_0)}\right)
\\ \displaystyle
=
\epsilon \sum_{m=1}^{n}(i_1-i_0) \varphi\left(\frac{m (u_{i_1}-u_{i_0})}{\epsilon(i_1-i_0)}\right)
-
\epsilon \sum_{m=1}^{n} (m-1) \varphi\left(\frac{m (u_{i_1}-u_{i_0})}{\epsilon(i_1-i_0)}\right).
\end{array}
\label{eq:QCP-energy-2}
\end{equation}
Here we see that the first term in the final expression for $E_2^{\rm QCP}$ is exactly the expression of the Cauchy-Born rule:
$\sum_{m=1}^{n} \varphi\left(\frac{m (u_{i_1}-u_{i_0})}{\epsilon(i_1-i_0)}\right)$ is the strain energy density associated with a representative atom and $\epsilon(i_1-i_0)$ is the length of the strained material.
We can also see that in the case of nearest neighbor interaction ($n=1$) the Cauchy-Born rule gives exactly the local component of the QCP energy and hence no corrections are required.
However, if $n>1$ then one needs to add the third component of energy (in order to avoid the ghost force), namely the energy of interaction of local and nonlocal atoms (atoms in the neighborhoods of $i_1$ and $i_K$):
\begin{equation}
\begin{array}{l} \displaystyle
E_3^{\rm QCP} =
\left.
\epsilon \sum_{\substack{j<i_1<i \\ 0<|i-j|\le n}} \varphi\left(\frac{u_i-u_j}{\epsilon}\right)\right|_{u_j = \frac{i_1-j}{i_1-i_0} u_{i_0} + \frac{j-i_0}{i_1-i_0}u_{i_1}}
\\ \displaystyle
+
\left.
\epsilon \sum_{\substack{i<i_K<j \\ 0<|i-j|\le n}} \varphi\left(\frac{u_i-u_j}{\epsilon}\right)
\right|_{u_j = \frac{i_{K+1}-j}{i_{K+1}-i_K} u_{i_K} + \frac{j-i_K}{i_{K+1}-i_K}u_{i_{K+1}}}
\\[1em] \displaystyle
= \epsilon \sum_{\substack{j<i_1<i \\ 0<|i-j|\le n}}
\varphi\left(\frac{u_i-u_{i_1}}{\epsilon} + \frac{(i_1-j) (u_{i_1}-u_{i_0})}{\epsilon(i_1-i_0)}\right)
\\ \displaystyle
+
\epsilon \sum_{\substack{i<i_K<j \\ 0<|i-j|\le n}}
\varphi\left(\frac{u_{i_K}-u_i}{\epsilon} + \frac{(j-i_K) (u_{i_{K+1}}-u_{i_K})}{\epsilon(i_{K+1}-i_{K})}\right).
\end{array}
\label{eq:QCP-energy-3}
\end{equation}
Here we denoted $i_{K+1} = i_1 + N$ and implied $u_{i_{K+1}} = u_{i_1}+1$.

\paragraph{finding solution}
Thus, we have written the interaction energy in terms of our unknowns $u_{i_1},\ldots,u_{i_K}$.
Adding it up with the external potential energy will yield the total potential energy of the QCP method $\Pi^{\rm QCP}$.
The QC solution is defined as a critical point $\Pi^{\rm QCP}$.
Finding a critical point typically consists in starting with some initial guess and performing Newton's iterations \eqref{eq:general_QCP_Newton}.
For Newton's iterations one needs the right-hand side vector $\left[\frac{\partial \Pi^{\rm QCP}}{\partial u_{i_k}}\right]_k$ and the stiffness matrix $\left[\frac{\partial^2 \Pi^{\rm QCP}}{\partial u_{i_k}\partial u_{i_l}}\right]_{kl}$, both of which can be computed from the above formulae in a rather straightforward manner.
We discuss the implementation of the projection method in subsection \ref{sec:QCP-implementation}.

%\paragraph{computing the force and the stiffness matrix}
%When computing the stiffness matrix, one can employ the FEM-type technique of assembing the matrix: instead of doing it entry-by-entry, one can do it atom pair-by-atom pair.
%More precisely, the suggested algorithm is as follows.
%\begin{enumerate}
%\item For each continuum element (in our case, only for element $i_0\le i\le i_1$), compute the interaction energies between atoms $i$
%\end{enumerate}

\paragraph{interpretation of quasinonlocal atoms}
First, note that in the case of nearest neighbor interaction ($n=1$), $E_2^{\rm QCP}$ is exactly the continuum energy computed by the Cauchy-Born rule and $E_3^{\rm QCP} \equiv 0$.
It means that the QCP method is equivalent to the traditional QC method in this case.

\paragraph{interpretation of quasinonlocal atoms}
Second, note that the local component of the energy $E_2^{\rm QCP}$ does not include interaction of nonlocal atoms with certain local atoms close to the interface (namely, local atoms $i_1-n+1,\ldots,i_1-1$ and $i_K+1,\ldots,i_K+n-1$).
These atoms play the role similar to the quasinonlocal atoms in the QNL method \cite{Shimokawa2004}, with the difference that in the original QNL method only nonlocal atoms change the way they interact, whereas in the QCP method the local atoms close to the interface change the way they interact.
The QNL method successfully removes the ghost force for second nearest neighbor interaction ($n=2$), but still exhibits the ghost force for longer interaction distances (though the magnitude of the ghost force of the QNL method is less than that of the QCE method, as we will see from the numerical results presented in section \ref{sec:results}).
%Below we will give a more detailed comparison of the QCP with the QNL method, as wellwill for second nearest neighbor interaction.

\subsection{Comparison of QCP with GCR}

\paragraph{intro}
In this subsection we will show that, in relation to the local-nonlocal interface, QCP is not equivalent to the particular instance of GCR that was presented in \cite{E2006}, but one can find an instance of GCR which would be equivalent to QCP.
However, GCR and QCP are essentially different in the local region in the way they treat interfaces of elements in the local region: GCR (as well as QNL) uses the Cauchy-Born extrapolation as in the standard QC, while QCP uses the projection.
The Cauchy-Born extrapolation and the projection give different results in the case of a non-equidistant lattice in the local region, since they treat interface between elements in a different way.
As was shown in \cite{linsinum}, the Cauchy-Born extrapolation introduces an additional error at element interfaces in the local region.
This additional error may be crucial, as we will see in subsection \ref{sec:2dtest}.

\paragraph{comparison for $n=2$}
One characteristic feature of the QCP method is that the QCP atomistic interaction energy coincides with the exact atomistic energy for the atomistic lattice satisfying the QC constraints \eqref{eq:QC_constraints}.
It therefore is not equivalent to the QNL method for $n=2$ (see the calculation in section 3.1 which shows the QNL energy is different from the exact atomistic energy) and therefore not equivalent to QNL for any $n\ge 2$, but is equivalent to GCR for $n=2$ as mentioned at the end of section 3.

\paragraph{comparison of QCP with GCR for $n=3$}
For the case of $n=3$, substitution of the values of coefficients for GCR from table I of \cite{E2006} and direct computation (which is straightforward, but too bulky to be presented here) yields
\[
\Etot^{\rm GCR} - \Etot^{\rm QCP} =
\frac{\epsilon}{2}\varphi\left(\frac{2 (u_{i_1+1} - u_{i_1-2})}{3 \epsilon} + \frac{u_{i_1+1}-u_{i_1}}{3 \epsilon}\right)
\]
\[
+
\frac{\epsilon}{2}\varphi\left(\frac{2 (u_{i_1-2} - u_{i_1+1})}{3 \epsilon} + \frac{u_{i_1-2}-u_{i_1-1}}{3 \epsilon}\right)
\]
\[
- \epsilon\varphi\left(\frac{u_{i_1+1}-u_{i_1-2}}{\epsilon}\right)
+ \textnormal{ respective terms for atoms around $i_K$}.
\]
As before, defining $D^2u_i = \frac{u_{i+1}-2 u_i+u_{i-1}}{\epsilon^2}$ and expanding $\Etot^{\rm QCP} - \Etot^{\rm GCR}$ in Taylor series w.r.t.\ $u_{i_2}$ around its extrapolated value $2 u_{i_1} - u_{i_1-1}$ (the respective terms for atoms around $i_K$ are expanded w.r.t.\ $u_{i_{K-1}}$) yields
\[
\Etot^{\rm QCP} - \Etot^{\rm GCR} =
%\frac{2\epsilon}{9} \varphi''\left(\frac{3 (u_{i_1}-u_{i_1-1})}{\epsilon}\right)
%\left(\frac{u_{i_2}-(1+i_1-i_0) u_{i_1} + u_{i_0}}{\epsilon(i_1-i_0)}\right)^2
%+
%\]
%\[
%O\left(\frac{u_{i_2}-(1+i_1-i_0) u_{i_1} + u_{i_0}}{\epsilon(i_1-i_0)}\right)^3
%=
%\]
%\[
%\frac{2\epsilon}{9} \varphi''\left(\frac{3 (u_{i_1}-u_{i_1-1})}{\epsilon}\right)
%\left(\frac{u_{i_2}-2 u_{i_1} + u_{i_1-1}}{\epsilon}\right)^2
%+
%O\left(\frac{u_{i_2}-2 u_{i_1} + u_{i_1-1}}{\epsilon}\right)^3
%=
%\]
%\[
\frac{2\epsilon}{9} \varphi''\left(\frac{3 (u_{i_1}-u_{i_1-1})}{\epsilon}\right)
\left(\epsilon D^2u_{i_1}\right)^2
+
O\left(\epsilon D^2u_{i_1}\right)^3
\]
\[
+ \frac{2\epsilon}{9} \varphi''\left(\frac{3 (u_{i_K}-u_{i_K+1})}{\epsilon}\right)
\left(\epsilon D^2u_{i_K}\right)^2
+
O\left(\epsilon D^2u_{i_K}\right)^3.
\]
One can check that the atomic force $\frac{\partial \Etot^{\rm GCR}}{\partial u_{i_k}}$ equals to $\frac{\partial \Etot^{\rm QCP}}{\partial u_{i_k}}$ up to the first order of accuracy in $\epsilon$ provided that the underlying solution is smooth (i.e.\ that $D^2 u_i$ is bounded).
For instance
\begin{equation}
\frac{\partial\Etot^{\rm QCP}}{\partial u_{i_2}} - \frac{\partial\Etot^{\rm GCR}}{\partial u_{i_2}} =
\frac{2}{9} \varphi''\left(\frac{3 (u_{i_1}-u_{i_1-1})}{\epsilon}\right)
O\left(\epsilon D^2u_{i_1}\right).
\label{eq:comparison_QCP_GCR_example}
\end{equation}
Moreover, the quantity
$\frac{2}{9} \varphi''\left(\frac{3 (u_{i_1}-u_{i_1-1})}{\epsilon}\right)$ is rather small in practice:
for instance in the case of Lennard-Jones potential, if the interatomic distance is approximately equal to $\epsilon$, then
\[
\frac{2}{9} \varphi''\left(\frac{3 (u_{i_1}-u_{i_0})}{\epsilon(i_1-i_0)}\right) \approx \frac{2}{9} \varphi''(3) \approx 0.003.
\]
Therefore we should expect that the difference in solutions by the QCP method and the GCR method is small, provided that the solution is smooth around the interface between the nonlocal and the local region (for $D^2 u_i$ not to be too large).
However, it should be noted here that if the solution is not smooth enough in the local region, then the QCP method may give significantly better results as compared to the GCR method (since they treat interfaces of finite elements in the local region in a different way).
This will be shown in the section with numerical tests.

\paragraph{there exists GCR which is equivalent to QCP}
Having established that the instance of the GCR method presented in table I of \cite{E2006} is not equivalent to the QCP method, it is nevertheless possible to come up with the GCR method which would be exactly equivalent to the QCP method.
For that, we should obtain the GCR method with the nonlocal atoms reconstructed as fully nonlocal atoms.
That is, the corresponding $C^{\rm GCR}_{ij}$ in \eqref{eq:GCR-reconstruction} should all be equal to one if $i$ is the nonlocal atom, except for the interface atom ($i_1$ or $i_K$ in our case) which is allowed to be less than one if $j$ is the local atom.
This can easily be attained by shifting the rows of table I to right by $n-1$ for $n>2$ (i.e.\ by assigning new values of $C^{\pm}_i$ to the old values of $C^{\pm}_{i-(n-1)}$ in notations of \cite{E2006}; note that for $n\le 2$ the method is already equivalent to QCP).

\paragraph{computations for $n=3$}
Computations below verify that for $n=3$ the following GCR method (obtained by shifting coefficients $C^{\rm GCR}_{ij}$ as described above)
\begin{equation}
\tilde{C}^{\rm GCR}_{ij} =
\left\{
\begin{array}{lll}
1 & & \textnormal{$(i,j)$ = $(i_1-1,i_1)$ or $(i_K+1,i_K)$}, \\
1 & & \textnormal{$(i,j)$ = $(i_1-1,i_1+1)$ or $(i_K+1,i_K-1)$}, \\
1 & & \textnormal{$(i,j)$ = $(i_1-2,i_1+1)$ or $(i_K+2,i_K-1)$}, \\
1 & & \textnormal{$(i,j)$ = $(i_1-1,i_1+2)$ or $(i_K+1,i_K-2)$}, \\
2/3 & & \textnormal{$(i,j)$ = $(i_1-3,i_1)$ or $(i_K+3,i_K)$}, \\
1/3 & & \textnormal{$(i,j)$ = $(i_1,i_1-3)$ or $(i_K,i_K+3)$}, \\
C^{\rm QCE}_{ij} & & \textnormal{otherwise}
\end{array}
\right.
\label{eq:GCR-tilde-coefficients}
\end{equation}
is equivalent to QCP.
Indeed, the internal energy for the QCP method for $n=3$ is
\[
\begin{array}{rcl}
\Etot^{\rm QCP} = \Etot^{\rm QCE}
& + &\displaystyle \left(
	\frac{\epsilon}{2} \varphi\left( \frac{u_{i_1-1}-u_{i_1+1}}{\epsilon} \right)
	-
	\frac{\epsilon}{2} \varphi\left( \frac{2 (u_{i_1-1}-u_{i_1})}{\epsilon} \right)
\right)
\\[1em]
& + &\displaystyle \left(
	\frac{\epsilon}{2} \varphi\left( \frac{u_{i_1-2}-u_{i_1+1}}{\epsilon} \right)
	-
	\frac{\epsilon}{2} \varphi\left( \frac{3 (u_{i_1-2}-u_{i_1-1})}{\epsilon} \right)
\right)
\\[1em]
& + &\displaystyle \left(
	\frac{\epsilon}{2} \varphi\left( \frac{u_{i_1-1}-u_{i_1+2}}{\epsilon} \right)
	-
	\frac{\epsilon}{2} \varphi\left( \frac{3 (u_{i_1-1}-u_{i_1})}{\epsilon} \right)
\right)
\\[1em]
& + &\displaystyle \textnormal{respective corrections for atoms around $i_K$}.
\end{array}
\]
The formula for the energy of this GCR method, according to \eqref{eq:GCR-tilde-coefficients} is then
\begin{equation}
\begin{array}{rcl}
\tilde{E}_{\rm tot}^{\rm GCR} = \Etot^{\rm QCE}
& + &\displaystyle \left(
	\frac{\epsilon}{2} \varphi\left( \frac{u_{i_1-1}-u_{i_1}}{\epsilon} \right)
	-
	\frac{\epsilon}{2} \varphi\left( \frac{u_{i_1-1}-u_{i_1}}{\epsilon} \right)
\right)
\\[1em]
& + &\displaystyle \left(
	\frac{\epsilon}{2} \varphi\left( \frac{u_{i_1-1}-u_{i_1+1}}{\epsilon} \right)
	-
	\frac{\epsilon}{2} \varphi\left( \frac{2 (u_{i_1-1}-u_{i_1})}{\epsilon} \right)
\right)
\\[1em]
& + &\displaystyle \left(
	\frac{\epsilon}{2} \varphi\left( \frac{u_{i_1-2}-u_{i_1+1}}{\epsilon} \right)
	-
	\frac{\epsilon}{2} \varphi\left( \frac{3 (u_{i_1-2}-u_{i_1-1})}{\epsilon} \right)
\right)
\\[1em]
& + &\displaystyle \left(
	\frac{\epsilon}{2} \varphi\left( \frac{u_{i_1-1}-u_{i_1+2}}{\epsilon} \right)
	-
	\frac{\epsilon}{2} \varphi\left( \frac{3 (u_{i_1-1}-u_{i_1})}{\epsilon} \right)
\right)
\\[1em]
& + &\displaystyle \frac{2}{3} \left(
	\frac{\epsilon}{2} \varphi\left( \frac{u_{i_1-3}-u_{i_1}}{\epsilon} \right)
	-
	\frac{\epsilon}{2} \varphi\left( \frac{3 (u_{i_1-3}-u_{i_1-2})}{\epsilon} \right)
\right)
\\[1em]
& + &\displaystyle \frac{1}{3} \left(
	\frac{\epsilon}{2} \varphi\left( \frac{u_{i_1}-u_{i_1-3}}{\epsilon} \right)
	-
	\frac{\epsilon}{2} \varphi\left( \frac{3 (u_{i_1}-u_{i_1-1})}{\epsilon} \right)
\right)
\\[1em]
& + &\displaystyle \textnormal{respective corrections for atoms around $i_K$}.
\end{array}
\label{eq:GCR-tilde-energy}
\end{equation}
Notice that the term in parenthesis in the first line of \eqref{eq:GCR-tilde-energy} is identically equals zero (this term was kept to match the formula \eqref{eq:GCR-tilde-coefficients}).
Taking into account that these energies are computed on an admissible QC lattice satisfying the constraints \eqref{eq:QC_constraints}, one can also notice that the terms in the fourth and the fifth line of \eqref{eq:GCR-tilde-energy} are identically zero, which makes $\tilde{E}_{\rm tot}^{\rm GCR} = \Etot^{\rm QCP}$ for $n=3$.
However, as we mentioned earlier, QCP and GCR are different if there are nodal atoms in the interior of the local region.
The subsection \ref{sec:2dtest} shows that this may be crucial in some cases.

\paragraph{important remark}
We would like to summarize the analysis of the QCP method with one important remark.
In estimates such as \eqref{eq:comparison_QCP_GCR_example}, following the conventions set in earlier works, we estimated the QC energy or the QC force through the non-nodal atoms such as $u_{i_1-1}$.
One should understand this in the following way.
The QC energy, as introduced in this work, depends only on positions of the nodal atoms $\uQC$.
Therefore, strictly speaking, the positions of the non-nodal atoms should be expressed through $\uQC$ using \eqref{eq:QC_constraints} and only then used in the estimates.
For instance, one should express $u_{i_1-1}$ as
\[
u_{i_1-1} = \frac{1}{i_1-i_0} u_{i_0} + \frac{i_1-i_0-1}{i_1-i_0} u_{i_1}.
\]
Thus, the estimate \eqref{eq:comparison_QCP_GCR_example} in its expanded form should read
\[
\frac{\partial\Etot^{\rm QCP}}{\partial u_{i_2}} - \frac{\partial\Etot^{\rm GCR}}{\partial u_{i_2}} =
\frac{2}{9} \varphi''\left(\frac{3 (u_{i_1}-u_{i_0})}{\epsilon (i_1-i_0)}\right)
O\left(\epsilon D^2u_{i_1}\right),
\]
where $D^2u_{i_1}$ should be expressed as
\[
D^2u_{i_1} = \frac{1}{\epsilon} \left(
	\frac{u_{i_2}-u_{i_1}}{\epsilon} - \frac{u_{i_1}-u_{i_0}}{\epsilon (i_1-i_0)}
\right).
\]

\subsection{On Implementation of the Projection Method}\label{sec:QCP-implementation}

\paragraph{intro}
This subsection is devoted to the implementation of the projection method.
We first describe the possible implementation in the 1D periodic setting and then discuss the implementation in the general case.
The implementation of the projection method consists in computing the right-hand side vector $\left[\frac{\partial \Pi^{\rm QCP}}{\partial u_{i_k}}\right]_k$ and the stiffness matrix $\left[\frac{\partial^2 \Pi^{\rm QCP}}{\partial u_{i_k}\partial u_{i_l}}\right]_{kl}$, which are the first and the second derivatives of the total energy with respect to the positions of the nodal atoms $u_{i_k}$.
The total energy of the projection method is $\Pi^{\rm QCP} = \Etot^{\rm QCP} + \Eext^{\rm QCP}$, where $\Etot^{\rm QCP}$ is defined by equations \eqref{eq:QCP-energy}--\eqref{eq:QCP-energy-3} and $\Eext^{\rm QCP}=\Eext$ is given by \eqref{eq:Eext}.
Computing derivatives of $\Eext^{\rm QCP}$ is fairly straightforward since it is a linear functional of our unknowns $u_{i_k}$.
That would be somewhat less straightforward if we assumed a smooth nonzero external force in the local region: in this case we need to express $\Eext$ through the positions of the nodal atoms $u_{i_k}$ and possibly use the approximate nodal summation techniques to optimize the performance of the method.

\paragraph{derivatives of $E_1^{\rm QCP}$, $E_2^{\rm QCP}$, and $E_3^{\rm QCP}$}
Differentiation of atomistic energy $E_1^{\rm QCP}$ is relatively straightforward: one should go through all atoms that contribute a nonzero interaction potential with the given atom and compute the respective entries of the matrix and the right-hand side vector.
Computation of derivatives of $E_2^{\rm QCP}$ is also staightforward: it is the standard QC energy with some correction terms.
The energy $E_2^{\rm QCP}$ contributes only to the entries associated with the nodal atoms in the local region ($i_1$ and $i_K$ in our case).
Computation of derivatives of $E_3^{\rm QCP}$ is somewhat less straightforward since $E_3^{\rm QCP}$ involves the interaction energy of the nodal and the non-nodal atoms.

\paragraph{computing derivatives of $E_3^{\rm QCP}$}
The strategy of computing derivatives of $E_3^{\rm QCP}$ is the following: we go through all pairs of atoms $u_{\alpha}$, $u_{\beta}$ that give a nonzero contribution to the interaction energy, where the atom $u_{\alpha}=u_{i_k}$ is in the nonlocal region and the atom $u_{\beta}$ is in the local region.
For the atom in the local region, we should represent its coordinate through our unknowns $u_{i_l}$:
\[
u_{\beta} = \boldsymbol{U}^T_{\beta} \left[u_{i_l}\right]_l,
\]
where $\left[u_{i_l}\right]_l$ is the vector of unknowns and $\boldsymbol{U}_{\beta}$ is the vector of coefficients of representation of $u_{\beta}$.
Note that the vectors $\boldsymbol{U}_{\beta}$ are essentially the columns of the matrix $U^{\rm QC}$ \eqref{eq:matrix-U}.
In the 1D periodic setting there are only two nonzero components of $\boldsymbol{U}_{\beta}$.
For atoms $\alpha$ and $\beta$ the term $\epsilon\varphi\left((u_{\alpha}-u_{\beta})/\epsilon\right)$ gives the following contribution to the right-hand side vector:
\[
\varphi'((u_{\alpha}-u_{\beta})/\epsilon)
\left(\boldsymbol{e}_k
-
\boldsymbol{U}_{\beta}
\right),
\]
and the following contribution to the stiffness matrix:
\[
\frac{1}{\epsilon} \varphi''((u_{\alpha}-u_{\beta})/\epsilon)
\left(\boldsymbol{e}_k \boldsymbol{e}^T_k
-
\boldsymbol{e}_k \boldsymbol{U}^T_{\beta}
-
\boldsymbol{U}_{\beta} \boldsymbol{e}^T_k
-
\boldsymbol{U}_{\beta} \boldsymbol{U}^T_{\beta}
\right)
\]
\[
=
\frac{1}{\epsilon} \varphi''((u_{\alpha}-u_{\beta})/\epsilon)
\left(\boldsymbol{e}_k - \boldsymbol{U}_{\beta} \right)
\left(\boldsymbol{e}_k - \boldsymbol{U}_{\beta} \right)^T.
\]
Here $\boldsymbol{e}_k$ is the unit vector corresponding to $u_{i_k}$.
Adding contributions of $E_3^{\rm QCP}$ to the previously computed contributions of $E_1^{\rm QCP}$ and $E_2^{\rm QCP}$ concludes assembling of the stiffness matrix and the right-hand side vector.

\paragraph{generalizations}
The described implementation methodology of the QCP method is not confined to the periodic 1D setting, but can easily be generalized to 2D or 3D, triangles whose sides are not aligned with the atomistic lattice, non-triangular mesh, etc.
For that, one just needs to go through all pairs of atoms and sum their contributions in the following way.
Consider two atoms $\alpha$ and $\beta$, whose positions are represented as
\[
u_{\alpha} = \boldsymbol{U}^T_{\alpha} \left[u_{i_k}\right]_k
\quad
\textnormal{and}
\quad
u_{\beta} = \boldsymbol{U}^T_{\beta} \left[u_{i_l}\right]_l.
\]
Then the contribution of the term $\epsilon\varphi\left((u_{\alpha}-u_{\beta})/\epsilon\right)$ to the stiffness matrix and the right-hand side vector can be easily expressed through the difference of their representation vectors $\boldsymbol{U}_{\alpha}-\boldsymbol{U}_{\beta}$.
One should take into account that the vectors $\boldsymbol{U}_{\alpha}$ and $\boldsymbol{U}_{\beta}$ are sparse and implement the respective vector and matrix operations efficiently.

\paragraph{cluster summation}
However, going through all the atoms in the summation is expensive and should be avoided.
In order to optimize the summation procedure, the following observation is to be made.
%Another observation which is crucial for optimization is to be made.
For all atoms $\alpha$ and $\beta$ within one element, all contributions to the right-hand side vector and the stiffness matrix can be computed once per element.
It is commonly implemented by choosing a representative atom well inside the element and computing its energy contribution\footnote{In a number of works the representative atoms where chosen at the element vertices, which has been shown recently to yield inaccurate results, see \cite{LuskinOrtner08} for further discussion.}.
Similar optimization can be done if the element boundaries (in 2D or 3D) are aligned with the atomistic lattice.
Then one can identify a fixed number (which depends on the interaction distance and the lattice structure) of representative atoms at element edges and faces.
Through these representative atoms all contributions given by atoms in different elements, in principle, can be computed once per mesh edge or face.
Thus, in such case the number of operations would only depend on the number of elements $K$ and would not depend on the actual number of atoms $N$ which can be very large.
However, if the element boundaries in a 2D or a 3D lattice are not aligned with the mesh, then the number of operations of the QCP method will be higher than $O(K)$.
In such cases one may consider employing some cluster summation rules which approximately sum all the contribution of atoms close to the boundary with $O(K)$ operations.

\subsection{QCP Method and Ghost Force}\label{sec:QCP-and-ghost-force}

\paragraph{intro}
In this subsection we show that the QCP method does not exhibit any ghost force whatsoever.
%The absence of the ghost force just follows from the fact that the uniform lattice ???
It should be noted that absence of the ghost force was first noticed by Rudd and Broughton in \cite{RB05} for the so-called coarse-grained molecular dynamics (CGMD) method in a general setting.
The QCP method is essentially a zero-temperature static CGMD method, therefore the argument of \cite{RB05} can be applied to it, with the modification that the argument of \cite{RB05} concerned the linearized model whereas we consider a general nonlinear model.
Here we adopt the Rudd and Broughton's argument to the QCP method in the 1D periodic setting as well as in the general setting.

\paragraph{setting for the ghost force}
Consider the zero external force $f_i = 0$.
Then obviously a uniform lattice $u_i = \epsilon z i$ satisfies the equilibrium equations \eqref{eq:original_problem}, which means that the total force on any atom is zero.
Let $i_k$  ($1\le k \le K$) define the nodal atoms for the QCP approximation.
Obviously a uniform lattice $u_i = \epsilon z i$ satisfies the QC constraints \eqref{eq:QC_constraints}, therefore it is an admissible QC lattice.
We say that a certain QC method does not exhibit the ghost force if $u_i = \epsilon z i$ is indeed the QC solution to the problem.
Otherwise we say that the QC approximation introduces some ghost forces that cause the uniform lattice to lose its equilibrium.
The following theorem states that QCP does not exhibit the ghost force.

\begin{theorem}[absence of ghost force for the 1D periodic setting]
Consider the zero external force $f_i = 0$ and a uniform lattice $u^0_i = \epsilon z i$.
Let $i_k$ ($1\le k \le K$) define the nodal atoms for the QCP approximation.
Then $\boldsymbol{u}^{0,{\rm QC}} = [u^0_{i_k}]_k$ is the QCP solution, i.e.
\begin{equation}
\left.\frac{\partial \Pi^{\rm QCP}}{\partial u_{i_k}}\right|_{u_{i_k}=u^0_{i_k}} = 0 \quad (1\le k \le K).
\label{eq:QCP-no-ghost-force}
\end{equation}
\end{theorem}

{\itshape Proof.}
First, notice that in the definition of $\Pi$ \eqref{eq:Pi-QCP} the term $U^{\rm QC}\boldsymbol{u}^{0,{\rm QC}}$ equals to $\boldsymbol{u}^{0}$, since the linear interpolation operator $U^{\rm QC}$ on the uniform lattice $\boldsymbol{u}^{0}$ is identical.
Second, in the 1D periodic setting the uniform lattice $\boldsymbol{u}^{0}$ is in equilibrium:
\[
\left.\frac{\partial \Pi}{\partial u_i}\right|_{u_i=u^0_i} = 0.
\]
Combining these two observations yields
\begin{equation}
\left.\frac{\partial \Pi^{\rm QCP}}{\partial u_{i_k}}\right|_{u_{i_k}=u^0_{i_k}} =
\frac{\partial \Pi(U^{\rm QC}\boldsymbol{u}^{0,{\rm QC}})}{\partial u_{i_k}} =
\frac{\partial \Pi(\boldsymbol{u}^{0})}{\partial u_{i_k}} = 0,
\label{eq:QCP-has-no-ghost-force}
\end{equation}
which concludes the proof of the theorem.

\paragraph{transition to the 2nd theorem}
We can observe that in the proof we essentially used only the fact that the uniform lattice is in equilibrium and at the same time is identity under the interpolation operator $U^{\rm QC}$.
This leads us to the following theorem in the general setting.

\begin{theorem}[absence of ghost force for the general setting] \label{thm:general-equilibrium}
Let $\Pi = \Etot + \Eext$ be the potential energy of a general atomistic system described by its degrees of freedom $\boldsymbol{u}\in \mathbb{R}^N$.
Let the lattice $\boldsymbol{u}^0$ be uniform in the sense that it is in equilibrium w.r.t.\ $\Pi$ (in other words, $\boldsymbol{u}^0$ is a critical point of $\Pi$).
Let the QC approximation $U^{\rm QC}$ be chosen so that there exists a QC vector $\boldsymbol{u}^{0,{\rm QC}}$ corresponding to the uniform lattice $\boldsymbol{u}^0 = U^{\rm QC} \boldsymbol{u}^{0,{\rm QC}}$.
Then $\boldsymbol{u}^{0,{\rm QC}}$ is a QC solution.
\end{theorem}

\paragraph{about the proof}
The proof of this theorem involves straightforward use of theorem conditions to show that \eqref{eq:QCP-has-no-ghost-force} also holds in the general setting.

\paragraph{stronger results}
In fact, a result somewhat stronger than Theorem \ref{thm:general-equilibrium} can easily be proven: if in addition $\boldsymbol{u}^0$ is a local minimizer of $\Pi$, then $\boldsymbol{u}^{0,{\rm QC}}$ is a local minimizer of $\Pi^{\rm QCP}$.
This relies on a simple fact that if a functional is restricted to a subspace which contains a local minimum $\boldsymbol{u}^0$, then $\boldsymbol{u}^0$ remains a local minimum of the restricted functional.

\section{Numerical Results}\label{sec:results}

\paragraph{intro}
First, we present the results of 1D tests in the same periodic setting under which we did analysis of the QCP method.
Then we present the results of a 2D test.
In all tests we chose the Lennard-Jones potential $\varphi(z) = z^{-12} - 2 z^{-6}$ with a cut-off radius of $3.25$.
Such potential involves third nearest neighbor interaction in one dimension.
The results by the four methods, namely: the projection method (QCP), the classical quasicontinuum method also known as the energy-based quasicontinuum method (QCE), the quasinonlocal quasicontinuum method (QNL), and the geometrically consistent reconstruction-based method (GCR) are presented and compared.
For the GCR method we took the particular instance of the method with coefficients given in table I of \cite{E2006} for the 1D tests and from table II for the 2D test.

\subsection{1D tests}

\subsubsection{Test with Localized External Force}

\paragraph{description of the test}
We considered a periodic lattice with a period of length $1$ with $N=10\,000$ atoms:
\[
u_{i+N} = u_i+1 \quad (-\infty<i<\infty).
\]
To model a defect of the atomistic material in 1D, we considered the external force to be zero for all atoms except atoms $N/2$ and $N/2+1$:
\[
f_i =
\left\{
\begin{array}{rcl}
-1 & & i=N/2 \\
1 & & i=N/2+1 \\
0 & & \textnormal{otherwise.}
\end{array}
\right.
\]
Correspondingly, the nonlocal region was chosen to be comprised of the following $m$ atoms (only even values of $m$ were chosen):
\[
\In = [ N/2-m/2+1,N/2+m/2 ].
\]
Since the external force is zero in the local region, we should expect the interatomic distances not to vary essentially in the whole local region and therefore no additional nodal atoms are required inside the local region.
Thus, the nodal atoms were chosen as
\[
i_1 = N/2-m/2+1, \quad i_2 = N/2-m/2+2, \quad \ldots, \quad i_m = N/2+m/2.
\]

\begin{figure}
\begin{center}
	\includegraphics{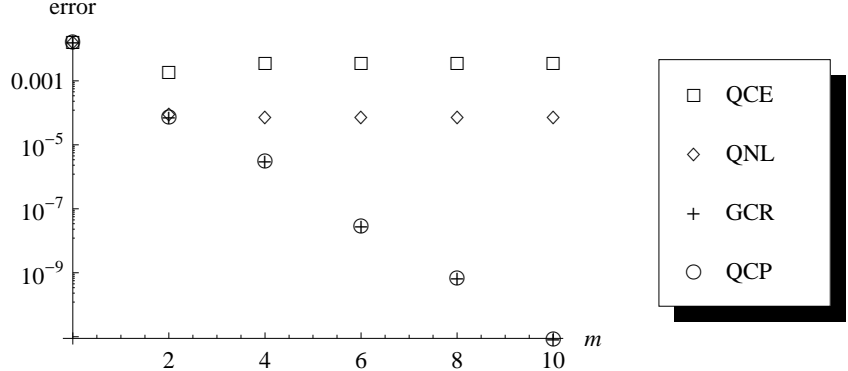}
\caption{Numerical error for 1D test with no bulk force. Error is computed in a discrete $W^{1,\infty}$ norm. $m$ is the number of nonlocal atoms near the defect.}
\label{fig:test1_error}
\end{center}
\end{figure}

\paragraph{results}
The dependence of the error in numerical solutions $\boldsymbol{u}^{\rm QCE}$, $\boldsymbol{u}^{\rm QNL}$, $\boldsymbol{u}^{\rm GCR}$, and $\boldsymbol{u}^{\rm QCP}$ on the length of the nonlocal region $m$ is presented in figure \ref{fig:test1_error}.
The error was computed as
\[
\textnormal{error} = |\boldsymbol{u}^{\rm numerical}-\boldsymbol{u}^{\rm exact}|_{W^{1,\infty}} =
\max\limits_{1\le i\le N} \left|\frac{(u_{i+1}^{\rm numerical}-u_{i+1}^{\rm exact})-(u_i^{\rm numerical}-u_i^{\rm exact})}{\epsilon}\right|.
\]
It can be seen that the errors of QCE and QNL do not converge to zero as $m$ grows, although the error of the QNL method is generally smaller.
Also, it can be seen that the errors of QCP and GCR converge to zero and are very close to each other, as expected from the above analysis.
The convergence rate is exponential in $m$.
The reason for the exponential convergence is that for the potential with a finite cut-off distance any disturbance of the uniform lattice decays exponentially (see, for example \cite{Dobson2009a-preprint} where a QC approximation is analyzed in the case of linearized problem with second nearest neighbor interaction).

\subsubsection{Test with Non-Localized External Force}

\paragraph{description of the test}
In the second test case we chose the same periodic lattice with $N=10\,000$, but the external force was taken to be the sum of the irregular and the regular component:
\[
\boldsymbol{f} = \boldsymbol{f}^{\rm irr} + \boldsymbol{f}^{\rm reg},
\]
where
\[
f^{\rm irr}_i =
\left\{
\begin{array}{rcl}
10 & & i=N/2 \\
-10 & & i=N/2+1 \\
0 & & \textnormal{otherwise,}
\end{array}
\right.
\]
and
\[
f^{\rm reg}_i = \frac{1}{N} \sin\left(1+\frac{2\pi i}{N}\right).
\]
The irregular component models the defect of material and the regular part models the external bulk force.
In this case the deformation in the local region is not constant and hence the nodal atoms should be chosen in the local region as well.
In the present test the nodal atoms were comprised of a number of equidistantly spaced atoms in the local region and all atoms in the nonlocal region.

\begin{figure}
\begin{center}
	\includegraphics{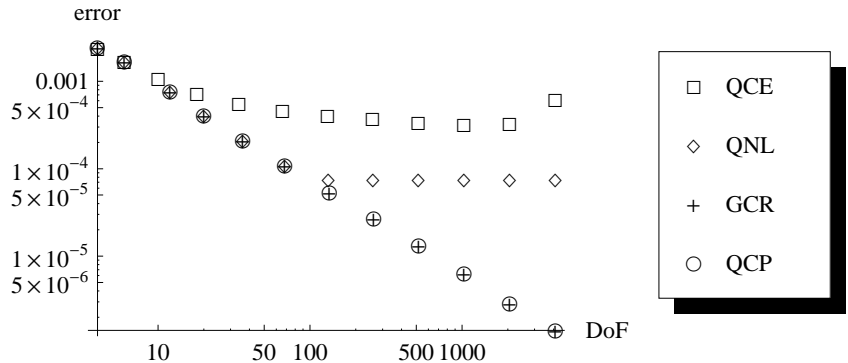}
\caption{Numerical error for 1D test with bulk force. Error is computed in a discrete $W^{1,\infty}$ norm. DoF is the number of degrees of freedom which equals to the total number of nodal atoms minus one.}
\label{fig:test2_error}
\end{center}
\end{figure}

\paragraph{results}
The dependence of errors in numerical solutions $\boldsymbol{u}^{\rm QCE}$,  $\boldsymbol{u}^{\rm QNL}$, $\boldsymbol{u}^{\rm GCR}$, and $\boldsymbol{u}^{\rm QCP}$ on the number of degrees of freedom (DoF) is presented in figure \ref{fig:test2_error}.
The number of degrees of freedom is the number of nodal atoms minus one (because we fix the position of one atom to avoid uncontrolled shift of the lattice as a whole).
For each value of DoF and each method we chose the length $m$ of the nonlocal region in such a way that it minimizes the error for the fixed DoF.
It can be seen that the QCP method and the GCR method again give very close results.
The errors of $\boldsymbol{u}^{\rm QCP}$ and $\boldsymbol{u}^{\rm GCR}$ are seen to decay with the first order in DoF.
The errors of QCE and QNL are seen to decay for small DoF and be approximately constant for large DoF because of dominating interface error.

\subsection{2D test}\label{sec:2dtest}

\begin{figure}
\begin{center}
\subfigure[full view]{
	\includegraphics{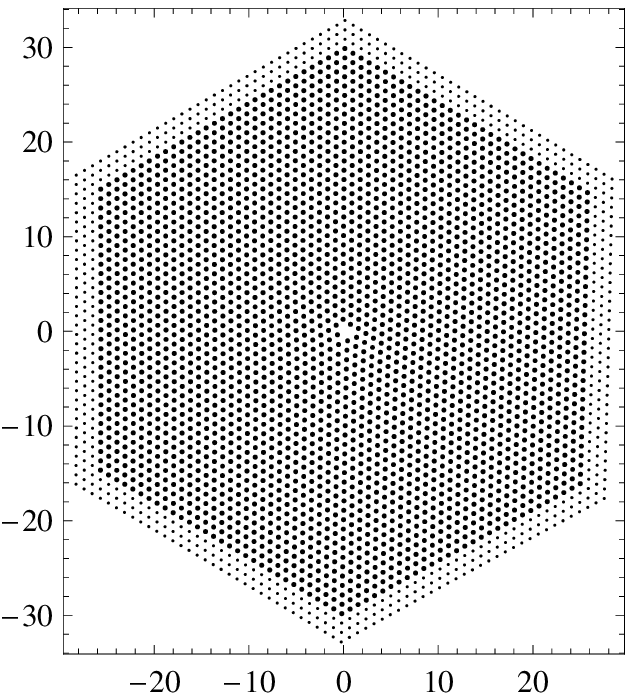}
	\label{fig:2dtest-configuration-full}
}
\hfill
\subfigure[zoomed view]{
	\includegraphics{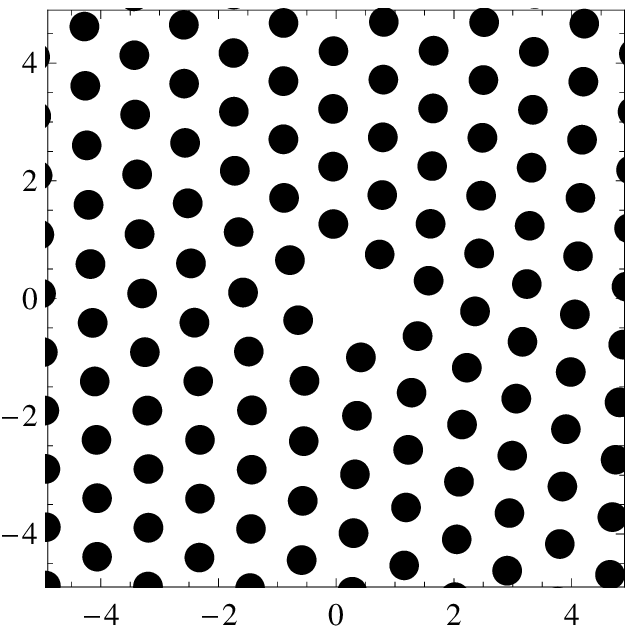}
	\label{fig:2dtest-configuration-zoom}
}
\end{center}
\caption{2D test with a point singularity, actual configuration}
\label{fig:2dtest-configuration}
\end{figure}

\begin{figure}
\begin{center}
\subfigure[full view]{
	\includegraphics{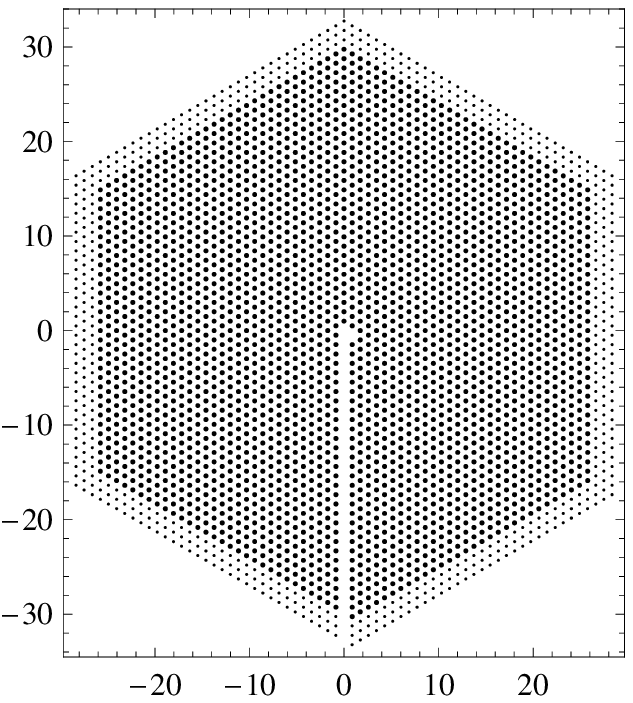}
	\label{fig:2dtest-reference-full}
}
\hfill
\subfigure[zoomed view]{
	\includegraphics{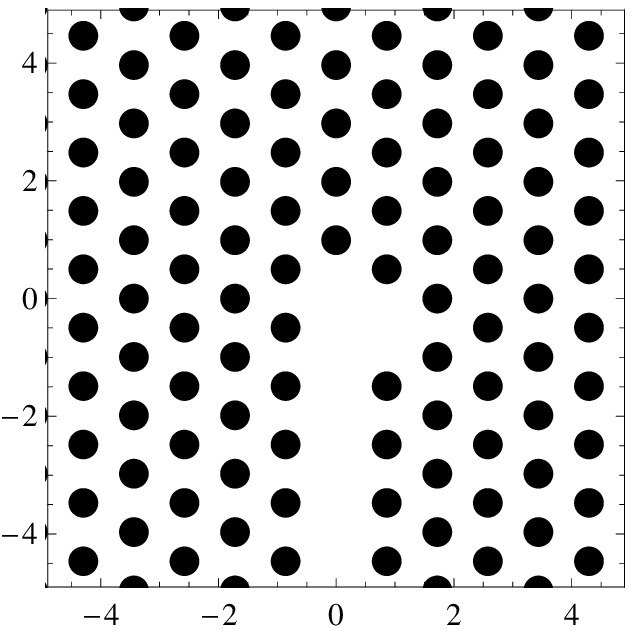}
	\label{fig:2dtest-reference-zoom}
}
\end{center}
\caption{2D test with a point singularity, reference configuration}
\label{fig:2dtest-reference}
\end{figure}

\begin{figure}
\begin{center}
	\includegraphics{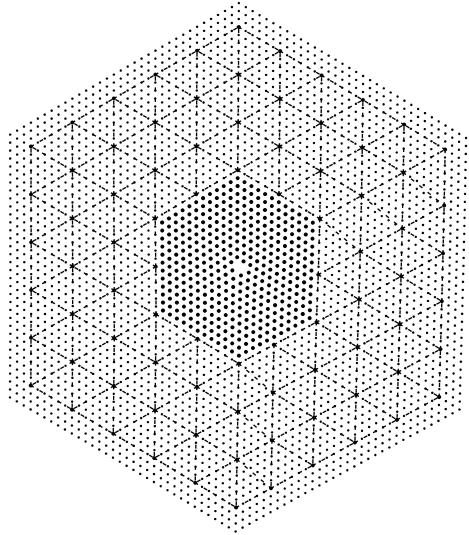}
\caption{2D test with a point singularity, a typical triangulation}
\label{fig:2dtest-triangulation}
\end{center}
\end{figure}

\paragraph{intro}
We now present a numerical test not under the one-dimensional periodic setting but under a general setting of the atomistic problem with defects.
We considered the configuration of atoms shown in fig.\ \ref{fig:2dtest-configuration} with the reference configuration shown in fig.\ \ref{fig:2dtest-reference}.
Such configuration models a point defect in a 2D lattice.
It corresponds to a dislocation in the respective 3D lattice (i.e.\ if one creates a 3D lattice by replicating these 2D layers, then the point defect will turn into a line defect which is an edge dislocation).
The lattice was comprised of $N_{\rm total}=3365$ atoms.
The behavior of the atomistic solution resembles the point singularity: the strain decays as $r^{-1}$ where $r$ is the distance to the defect (the $r^{-1}$ asymptotics of the strain is a general behavior of the strain for the edge dislocation \cite{Hull2001}.

\paragraph{boundary conditions}
The boundary conditions were set as follows.
First, the full atomistic problem with $N_{\rm total}=3365$ atoms was solved using the Neumann boundary conditions (i.e.\ the positions of the atoms on the boundary were not fixed).
This solution was used to fix the coordinates of the three outmost layers of atoms (shown as smaller circles in fig.\ \ref{fig:2dtest-configuration-full} and \ref{fig:2dtest-reference-full}).
Thus, the Dirichlet boundary conditions were used in the computations with three layers of atoms fixed.
This was done in order to eliminate the effect of the boundary conditions on the lattice near the boundary and thus to eliminate the effect of the free boundary on the error of the QC solution.
The typical triangulation used is shown in fig.\ \ref{fig:2dtest-triangulation}.
Notice that the three outmost layers of atoms do not participate in the triangulation.
In fact, they also do not participate in the atomistic solution, only yielding the external force on the boundary atoms.
Thus, the effective number of atoms in the computations was $N=2789$.

\begin{figure}
\begin{center}
	\includegraphics{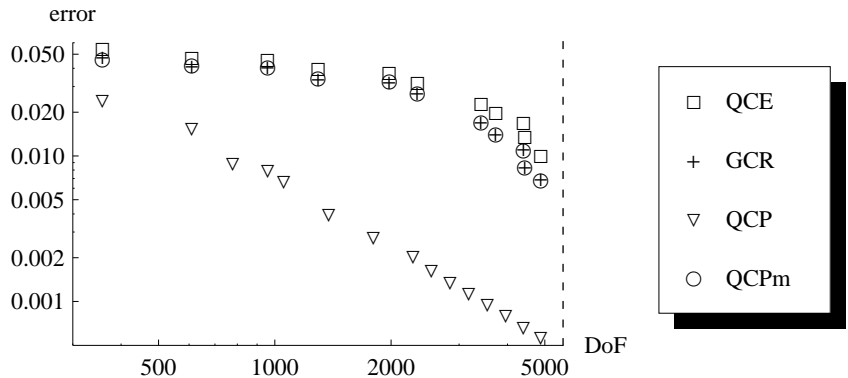}
\caption{Numerical error for the 2D test. Error is computed in an analog of $W^{1,\infty}$ norm. The dashed vertical line corresponds to the total number of degrees of freedom in the atomistic system.}
\label{fig:2dtest-error2}
\end{center}
\end{figure}

\paragraph{norm of error}
By similarity with the 1D tests, the error of numerical solution was computed in a vector $W^{1,\infty}$-norm, which is defined in the following way.
We loop over all triplets of neighboring atoms in the numerical and the exact solution and compute the Jacobian of the mapping of the triangle (corresponding to these three atoms) in the exact solution to the corresponding triangle in the numerical solution.
The 2-norm of the difference of this Jacobian with the unit matrix was taken as the error for one triangle.
Then we take the maximum over all such triplets of neighboring atoms, which is the sought $W^{1,\infty}$-norm.

\paragraph{results}
The dependence of the error in numerical solutions by different methods on number of degrees of freedom (DoF) is presented in figure \ref{fig:2dtest-error2}.
It can be clearly noticed that the errors of the QCP method are much less than those of the other methods and the error of the original QCE method is the highest.
We argue that the reason for the error of the QCP method to be much less than the errors of other methods is that it handles interfaces between the mesh triangles in the local region in a more accurate way.
To verify this, we also implemented the modified QCP method (QCPm) which is the same as QCP method everywhere except at the interfaces between the mesh triangles in the local region, where the standard QC discretization was used.
In other words, QCPm uses the QC discretization everywhere in the local region except at the interface between the nonlocal and local region, where QCPm behaves as QCP.
The results on fig.\ \ref{fig:2dtest-error2} show that the errors of the QCPm method and the GCR method are very close to each other.
It can also be seen that errors for GCR and QCPm remain constant for small and moderate DoF and start to decay only when the DoF of the numerical solution is close to DoF of the original atomistic problem (the total DoF of the atomistic problem is shown with the dashed line in fig.\ \ref{fig:2dtest-error2}).
On the contrary, the error of the QCP method decays steadily as DoF increases.

\section{Discussion and Conclusion}\label{sec:discussion-conclusion}

\paragraph{on different interaction distance}
The presented QCP method coincides with the classical QC approximation for nearest neighbor interaction ($n=1$).
It is not equivalent to the QNL method, although both methods effectively eliminate the ghost force for second nearest neighbor interaction ($n=2$).
As compared with the GCR method, we may conclude that QCP is different from the particular instance of the GCR method presented in \cite{E2006}, but there exist other instances of the GCR method which are equivalent to the QCP method (more precisely, to the QCPm method) in the way how these methods treat the local-nonlocal interface.

\paragraph{general trend of numerical error}
The results of 1D and 2D numerical tests show the following patterns.
The error of the QCP method was always lower than the errors of QCE and QNL.
The error of the QNL method was lower than that of the QCE method, although both errors do not converge to zero.

\paragraph{comparison with E, Lu, and Yang (2006)}
The error of the QCP method was essentially the same as the error of the GCR method for 1D tests and lower for the 2D tests.
We argue that the reason for such behavior is that in 1D case the disturbance of the uniform lattice decays exponentially, whereas the disturbance in the 2D test decays much slower: as the inverse of the distance to the defect.
This causes large interface errors between elements within the local region for all methods except QCP since the former use the Cauchy-Born extrapolation.
%, since the later does not directly use the Cauchy-Born extrapolation in the local region.
Not using the Cauchy-Born extrapolation in the local region results in a somewhat larger number of operations for assembling the stiffness matrix and the right-hand side vector, but may give much more accurate results as compared to the other methods.
In case if such increased accuracy is not essential then one can use the modified QCP (i.e.\ QCPm) method which treats the interface between nonlocal and local region as the QCP method and uses the Cauchy-Born extrapolation in the local region.
The QCPm method yields essentially the same results as the GCR method but, as we believe, is easier to implement, analyze, and generalize to other problems.

\paragraph{discussion on QCP and GCR}
GCR is a cleverly designed method and gives a very general criterion to remove the ghost force at the local-nonlocal interface.
However, based on the results of analysis and numerical experiments we infer that the QCP method (including QCPm, which is a simpler but less accurate modification of the projection method) is easier to implement and analyze, does not depend on any free parameters (such as reconstruction coefficients), and may be even more accurate than GCR in general,
in spite of the fact that QCP may be expressed as a particular case of GCR (with a difference that they treat the interface between elements in local regions differently).
First, as discussed in subsection \ref{sec:QCP-implementation}, for implementation of QCP one just needs the atom position's representations through the nodal atoms.
No other geometric information is needed for the implementation, which makes it easier than GCR, since the later requires the a priori tabulated coefficients of reconstruction and needs to determine distance to the local-nonlocal interface for each atom near the interface.
Second, the QCP method is also easy to analyze: one can benefit from the well-developed theory of finite elements which offers a powerful method of reducing the problem of convergence to the problem of approximation (see e.g.\ \cite{linsinum,Ortner2008}).
Third, the QCP method does not require determining any problem-dependent parameters beforehand.
It is therefore more flexible: one needs to do less investment for solving the new problems.
For example, if one wants to apply the QC method for metallic alloys, whose atomistic lattices are not uniform due to presence of atoms of different metals in the lattice, one just needs to describe this lattice in terms of reconstruction of the non-nodal atoms through the nodal atoms.
It should also be noted that the QCP method is essentially a particular case of the CGMD method originally proposed for coarse-graining the finite temperature molecular dynamics \cite{rudd2, RB2000, RB05}.
Therefore it may potentially have wider applications.
Fourth, as we have seen from the 2D numerical test (subsection \ref{sec:2dtest}), the QCP method may be more accurate since it also offers a natural way to eliminate the edge error in the nonlocal region as compared to the standard Cauchy-Born extrapolation.

\section*{Acknowledgements}
P.\ Lin would like to thank for the research support from Leverhulme Trust (Grant number RF/9/RFG/2009/0507).

\bibliographystyle{siam}
\bibliography{paper3}

\end{document}